\documentclass[11pt]{article}

\usepackage[margin=1in]{geometry}
\usepackage[T1]{fontenc}

\usepackage{amsmath}
\usepackage{amsfonts}
\usepackage{amsthm}
\usepackage{mathtools}

\usepackage{newtxtext}
\usepackage{newtxmath}

\usepackage{microtype}

\usepackage{graphicx}
\usepackage{caption}
\usepackage{subcaption}

\usepackage{algorithm}
\usepackage{algpseudocode}

\usepackage{booktabs}
\usepackage{array}

\usepackage{enumitem}
\usepackage{listings}
\usepackage[normalem]{ulem}
\usepackage{xcolor}

\usepackage[hidelinks]{hyperref}
\usepackage[nameinlink,noabbrev]{cleveref}

\theoremstyle{definition}

\theoremstyle{remark}

\usepackage[hidelinks]{hyperref}
\usepackage[nameinlink,noabbrev]{cleveref}

\title{\textbf{3D Uncertainty Quantification for Photoacoustic Tomography}}

\author{
Babak Maboudi Afkham\thanks{
Unit of Mathematical Sciences, Faculty of Science,
University of Oulu, Oulu, Finland.
Corresponding author: \texttt{babak.maboudi@oulu.fi}
}
\and
Amal Mohammed A. Alghamdi\thanks{
College of Computing and Mathematics,
King Fahd University of Petroleum and Minerals.
Former affiliation: Department of Applied Mathematics and Computer Science,
Technical University of Denmark, Kgs.\ Lyngby, 2800, Denmark.
}
\and
Hassan Yazdanian\thanks{
Unit of Mathematical Sciences, Faculty of Science,
University of Oulu, Oulu, Finland.
}
\and
Tanja Tarvainen\thanks{
Department of Technical Physics,
University of Eastern Finland, Kuopio, Finland.
}
}

\date{}

\begin{document}

\maketitle

\begin{abstract}
Photoacoustic tomography (PAT) is a hybrid imaging modality that has emerged as a promising tool for high-resolution biomedical imaging, motivating the need for reliable uncertainty quantification (UQ) to assess the confidence of reconstructed images. Bayesian approaches provide a rigorous framework for UQ but remain computationally challenging for realistic three-dimensional PAT and are sensitive to numerical approximations introduced by the discretization of the governing wave equation. In this work, we develop a finite-element Bayesian UQ framework for PAT that naturally accommodates complex computational domains and detector geometries while enabling large-scale three-dimensional inference. The proposed methodology reformulates the randomize-then-optimize (RTO) sampling strategy as a matrix-free algorithm that generates independent posterior samples using only forward and adjoint wave propagations. Particular attention is given to the construction of an adjoint discretization that forms an exact transpose pair with the discrete forward operator while remaining consistent with the continuous PAT adjoint, thereby enabling efficient least-squares solvers within the sampling procedure. We further investigate the influence of temporal discretization, artificial boundary conditions, and adjoint consistency on posterior uncertainty, identifying discretization strategies that avoid numerical artifacts in uncertainty estimates. The proposed framework is validated through comparisons with exact posterior statistics, existing Bayesian PAT methods, and Hamiltonian Monte Carlo using the No-U-Turn Sampler (NUTS), and is demonstrated on a three-dimensional problem with approximately $2\times10^5$ unknowns on a general finite-element domain. To the best of our knowledge, this is the first large-scale Bayesian PAT study on general three-dimensional finite-element geometries, and the resulting methodology extends naturally to a broad class of linear PDE-constrained inverse problems.
\end{abstract}


% ============================================================
\section{Introduction}

photoacoustic tomography (PAT) is a hybrid imaging technique used to generate high-resolution 3D images of biological tissues \cite{li2009,beard2011} This method relies on the photoacoustic effect, where short laser pulses are absorbed by tissues, causing localized heating and a subsequent increase in pressure. This pressure is then released as ultrasonic waves, which are detected by ultrasound sensors and used to reconstruct a 3D image of the tissue \cite{li2009}.

The increasing adoption of photoacoustic tomography (PAT) in preclinical and clinical imaging has made uncertainty quantification (UQ) an increasingly important challenge. Multiple sources of uncertainty can compromise the accuracy and reliability of reconstructed images, ultimately affecting their diagnostic utility. These include measurement noise, limited-view data acquisition resulting from practical constraints on ultrasound sensor placement, and incomplete knowledge of tissue optical and acoustic properties \cite{doi:10.1137/18M1231341,hauptmann2020deep,tarvainen2024quantitative}. These challenges have motivated substantial efforts to characterize and quantify uncertainty in PAT. Ren et al. \cite{doi:10.1137/18M1231341} conducted a mathematical analysis of uncertainty sources in quantitative PAT and performed sensitivity analyses to assess their effects on image reconstruction. Bayesian approaches to uncertainty quantification have also been developed for PAT and quantitative PAT \cite{tarvainen2012reconstructing,pulkkinen2014,tick2016,tick2019modelling,sahlstrom2020modeling}. Although these methods provide rigorous uncertainty estimates, they are often computationally expensive, particularly for the large-scale, high-resolution inverse problems encountered in PAT.

Recent work has explored learning-based approaches for efficient UQ in PAT, including variational autoencoders (VAEs) \cite{Goh2021,sahlstrom2023utilizing} and Monte Carlo dropout \cite{seoni2025exploring,hauptmann2020deep}. While these methods can substantially reduce computational cost, they provide data-driven uncertainty estimates whose reliability depends on the quality and representativeness of the training data and may be affected by training biases or distribution shifts. Moreover, Monte Carlo dropout provides only approximate uncertainty estimates and has shown limited performance for initial pressure reconstruction.

A standard model for acoustic propagation in photoacoustic tomography (PAT) is a linear hyperbolic wave equation \cite{xia2014photoacoustic} with the initial pressure distribution as an unknown initial condition and time-dependent pressure measurements recorded on an observation surface. Although Bayesian uncertainty quantification for PDE-based inverse problems is well established conceptually, its reliable and computationally tractable implementation remains challenging. First, evaluation of the forward map requires spatial and temporal discretization and, on truncated computational domains, an approximation of the radiation condition. Dispersion, dissipation, numerical diffusion, and boundary reflections introduced by these approximations can contaminate the inferred posterior; unless the numerical error is controlled or incorporated into the statistical model, the reported uncertainty may reflect discretization error in addition to physical and observational uncertainty. Second, repeated numerical solution of the wave equation can be computationally expensive. Markov chain Monte Carlo (MCMC) methods \cite{kaipio05} provide a general, asymptotically exact approach to posterior sampling, but typically require many sequential likelihood evaluations, each involving one or more forward PDE solves. Consequently, conventional MCMC may be computationally prohibitive for large-scale, high-resolution PAT problems. These observations suggest that practical Bayesian uncertainty quantification for PAT requires not only computationally efficient sampling algorithms but also numerical discretizations that faithfully represent the underlying continuous inverse problem.

In this work, we develop a finite-element Bayesian uncertainty quantification framework for PAT that naturally accommodates complex computational domains and detector geometries and investigate the influence of numerical approximations on posterior inference. In particular, we study the effects of temporal discretization, artificial boundary conditions, and discretization of the forward, adjoint, and reverse-time wave equations on posterior uncertainty. These investigations identify strategies that avoid numerical artifacts in posterior uncertainty estimates. To reduce the computational cost of Markov chain Monte Carlo methods, we reformulate the randomize-then-optimize (RTO) approach \cite{bardsley2014randomize} for large-scale PAT. For the linear-Gaussian model considered here, the proposed formulation generates independent posterior samples while remaining entirely matrix-free, requiring neither assembled system matrices nor direct linear solves. Each sample requires only forward and adjoint wave propagation and inexpensive diagonal operations, making the method suitable for large-scale and high-performance computing.

A central contribution of this work is the construction of an adjoint discretization tailored for Bayesian posterior sampling in PAT. In PDE-constrained inverse problems, the discretization of the continuous (physical) adjoint generally does not coincide with the algebraic adjoint of the discrete forward operator, leading to so-called \emph{unmatched adjoint} problems \cite{dong2019fixing,huber2025novel}. This distinction is particularly important for optimization-based sampling methods such as RTO, where each posterior sample requires the solution of a least-squares problem. Efficient Krylov methods such as LSQR \cite{paige1982lsqr} rely on the adjoint being the exact transpose of the discrete forward operator; otherwise, the underlying least-squares problem is no longer solved by LSQR, requiring more general methods such as GMRES \cite{hansen1998rank}, which accommodate unmatched operators but are generally less efficient for repeatedly solving the large least-squares problems arising during posterior sampling. We therefore construct forward and adjoint discretizations that form an exact transpose pair with respect to the discrete inner product while remaining consistent with the continuous PAT adjoint under mesh refinement, enabling efficient use of LSQR within the Bayesian sampling algorithm.

The proposed framework is demonstrated on a three-dimensional finite-element PAT problem with approximately $2\times10^5$ unknowns on a general computational domain with nontrivial measurement geometries. To the best of our knowledge, this is the first large-scale Bayesian PAT study based on finite-element discretization, extending Bayesian PAT beyond Fourier pseudospectral frameworks predominantly restricted to Cartesian or relatively simple domains. The methodology is validated against exact posterior variances for small-scale problems, existing uncertainty quantification methods for two-dimensional PAT, and Hamiltonian Monte Carlo using the No-U-Turn Sampler (NUTS). Despite the Gaussian posterior, these experiments show that a standard NUTS implementation is computationally impractical at this scale, highlighting the importance of exploiting problem structure in Bayesian sampling. Finally, we investigate the influence of prior modeling in three dimensions. Although developed for PAT, the methodology and many associated numerical considerations extend to a broad class of large-scale linear PDE-constrained inverse problems.

This paper is organized as follows. \Cref{sec:PA,sec:FEM} present the PAT forward and inverse problems and their finite-element discretization. \Cref{sec:Bayes,sec:RTO} introduce the Bayesian formulation and RTO posterior sampling. \Cref{sec:results} presents numerical validation and large-scale results, followed by conclusions in \Cref{sec:conclusion}.

\section{The Photo-Acoustic Tomography Inverse Problem} \label{sec:PA}
In this section, we introduce the mathematical model governing acoustic wave propagation in photoacoustic tomography and formulate the corresponding inverse problem.

In PAT, a short laser pulse illuminates the tissue, with the deposited optical energy regarded as instantaneous \cite{xia2014photoacoustic}, generating an initial pressure distribution through the photoacoustic effect. This pressure subsequently propagates as an acoustic wave and is recorded by ultrasound sensors at the boundary of the imaging domain.
Under the thermal and stress confinement assumptions \cite{cox2005fast}, the absorbed optical energy can therefore be modeled as an initial pressure distribution. The subsequent acoustic propagation is described by the initial-value problem
\begin{equation} \label{eq:wave}
    \begin{aligned}
        \frac{\partial^2 p(\boldsymbol x,t)}{\partial t^2}
        &= c^2 \Delta p(\boldsymbol x,t),
        \qquad && \boldsymbol x \in \Gamma,\quad t\in[0,T], \\
        p(\boldsymbol x,0)
        &= p_0(\boldsymbol x),
        && \boldsymbol x\in \Gamma,\\
        \frac{\partial p(\boldsymbol x,0)}{\partial t}
        &= 0,
        && \boldsymbol x\in \Gamma .
    \end{aligned}
\end{equation}
Here $\Gamma \subset \mathbb R^d$, with $d=2$ or $3$, denotes the bounded computational domain, $\boldsymbol x\in\mathbb R^d$ denotes the spatial coordinate, $t\in[0,T]$ denotes time, and $T>0$ is the final observation time. The function $p:\Gamma\times[0,T]\to\mathbb R$ denotes the acoustic pressure, $p_0:\Gamma\to\mathbb R$ is the initial pressure distribution, and $c>0$ is the acoustic wave speed. Throughout this work, we consider the non-dimensionalized form of \eqref{eq:wave}, in which the wave speed is normalized to $c=1$.

In practice, acoustic waves generated in PAT propagate in an effectively unbounded medium and satisfy a radiation condition at infinity \cite{engquist1979radiation}. Since simulations use a bounded computational domain, an artificial boundary condition is required to approximate outgoing-wave behaviour while minimizing artificial reflections. Here, we employ the first-order absorbing boundary condition
\begin{equation} \label{eq:boundary}
    \frac{\partial p(\boldsymbol x,t)}{\partial \mathbf n}
    =
    -\frac{1}{c}\frac{\partial p(\boldsymbol x,t)}{\partial t},
    \qquad
    \boldsymbol x\in\partial\Gamma,
\end{equation}
where $\mathbf n$ denotes the outward unit normal on $\partial\Gamma$, the boundary of $\Gamma$. This boundary condition allows outgoing waves to leave the computational domain with minimal reflection; this choice is further justified in \Cref{sec:implicit-transpose}. We denote the solution operator associated with \eqref{eq:wave}--\eqref{eq:boundary} by $\mathcal A : p_0 \mapsto p$, where $p$ solves \eqref{eq:wave}--\eqref{eq:boundary}.

The pressure field is observed by $N_{\mathrm{obs}}$ ultrasound sensors located at positions $\boldsymbol  x_{\mathrm s}^i\in\Gamma_{\mathrm{obs}}\subseteq\partial\Gamma$, for $i=1,\ldots,N_{\mathrm{obs}}$. We define the observation operator $\mathcal O:\Gamma\times(0,T]\rightarrow\mathbb R^{N_T\times N_{\mathrm{obs}}}$, which maps the pressure field to the discrete measurements
\begin{equation}
    \begin{aligned}
    \mathbf Y &:= \mathcal O(p),\\
    [\mathbf Y]_{j,i} &= p(\boldsymbol x_{\mathrm s}^i,t_j), \qquad t_j=j\Delta t,\quad j=1,\ldots,N_T,
    \end{aligned}
\end{equation}
where $[\mathbf Y]_{j,i}$ denotes the pressure measured by the $i$th sensor at time $t_j$, $\Delta t$ is the temporal sampling interval, and $N_T$ is the total number of time measurements.

Ultrasound measurements are inevitably contaminated by measurement noise. In this work, we assume an additive Gaussian noise model and measure the noise level relative to the discrete approximation of the $L^2(0,T;\ell^2)$ norm of the complete space--time data. Specifically, the observed data are given by $\mathbf Y_{\mathrm{obs}}=\mathbf Y+\eta$,where
\[
[\eta_0]_{j,i}\stackrel{\mathrm{i.i.d.}}{\sim}\mathcal N(0,1),
\qquad
j=1,\ldots,N_T,\quad
i=1,\ldots,N_{\mathrm{obs}},
\]
and
\begin{equation} \label{eq:noise-model}
    \eta
    =
    d_{\mathrm{noise}}
    \frac{\|\mathbf Y\|_{L^2(0,T;\ell^2)}}
         {\|\eta_0\|_{L^2(0,T;\ell^2)}}
    \eta_0, \qquad 
    \|\mathbf Y\|_{L^2(0,T;\ell^2)}
    :=
    \left(
        \Delta t
        \sum_{j=1}^{N_T}
        \sum_{i=1}^{N_{\mathrm{obs}}}
        [\mathbf Y]_{j,i}^2
    \right)^{1/2},
\end{equation}
where $d_{\mathrm{noise}}>0$ denotes the prescribed relative noise level. We define the forward operator $\mathcal G := \mathcal O\circ\mathcal A$,
which maps the initial pressure distribution to the corresponding space--time measurements. The complete observation model is then given by
\begin{equation} \label{eq:forward}
    \mathbf Y_{\mathrm{obs}}
    =
    \mathcal G(p_0)
    +
    \eta,
\end{equation}
where $\eta$ is the additive Gaussian noise defined above. Since both the solution operator $\mathcal A$ and the observation operator $\mathcal O$ are linear, the forward operator $\mathcal G$ is also linear.

\section{Finite Element Discretization} \label{sec:FEM} 
In this section, we derive the finite element discretization of the forward problem. For temporal discretization, we rewrite \eqref{eq:wave} as a first-order system by introducing the auxiliary velocity $v=\partial p/\partial t$. This formulation enables the second-order St\"ormer--Verlet time integration scheme \cite{hairer2006geometric}, yielding the initial-boundary value problem

\begin{equation} \label{eq:first-order}
    \begin{aligned}
        \frac{\partial v(\boldsymbol x,t)}{\partial t}
        - c^2\Delta p(\boldsymbol x,t)
        &=0,
        && \boldsymbol x\in\Gamma,\quad t\in[0,T],\\
        \frac{\partial p(\boldsymbol x,t)}{\partial t}
        -v(\boldsymbol x,t)
        &=0,
        && \boldsymbol x\in\Gamma,\quad t\in[0,T],\\
        p(\boldsymbol x,0)
        &=p_0(\boldsymbol x),
        && \boldsymbol x\in\Gamma,\\
        v(\boldsymbol x,0)
        &=0,
        && \boldsymbol x\in\Gamma,\\
        \frac{\partial p(\boldsymbol x,t)}{\partial\mathbf n}
        +\frac1c\,v(\boldsymbol x,t)
        &=0,
        && \boldsymbol x\in\partial\Gamma,\quad t\in[0,T].
    \end{aligned}
\end{equation}

Let $w,q\in H^1(\Gamma)$ be test functions. Multiplying \eqref{eq:first-order} by $w$ and $q$, integrating over $\Gamma$, applying Green's identity to the first equation, and incorporating boundary conditions yields the weak formulation
\begin{equation}\label{eq:weak-form}
\left\{
\begin{aligned}
    \int_\Gamma v_t\,w\,d\boldsymbol x
    + c^2\int_\Gamma\nabla p\cdot\nabla w\,d\boldsymbol x
    +\frac1c\int_{\partial\Gamma}vw\,ds
    &=0,\\
    \int_\Gamma p_t\,q\,d\boldsymbol x
    -\int_\Gamma v\,q\,d\boldsymbol x
    &=0.
\end{aligned}
\right.
\end{equation}
Here, the subscript $t$ denotes differentiation with respect to time.

Let $V_h\subset H^1(\Gamma)$ denote the space of continuous piecewise linear finite element functions with basis $\{\phi_i\}_{i=1}^{N_h}$. Seeking the approximations $p_h(\boldsymbol x,t)=\sum_{i=1}^{N_h}p_i(t)\phi_i(\boldsymbol x)$, and $v_h(\boldsymbol x,t)=\sum_{i=1}^{N_h}v_i(t)\phi_i(\boldsymbol x)$,
and choosing the test functions from the same space yields the semi-discrete system
\begin{equation}\label{eq:semi-discrete}
\begin{aligned}
    \mathbf M\frac{d\mathbf v}{dt}
    +\mathbf K\mathbf p
    +\mathbf B\mathbf v
    &=0,\\
    \mathbf M\frac{d\mathbf p}{dt}
    -\mathbf M\mathbf v
    &=0,
\end{aligned}
\end{equation}
where
\begin{equation}
[\mathbf M]_{mn}=\int_\Gamma \phi_m\phi_n\,d\boldsymbol x,\;
[\mathbf K]_{mn}=c^2\int_\Gamma \nabla\phi_m\cdot\nabla\phi_n\,d\boldsymbol x,\;
[\mathbf B]_{mn}=\frac{1}{c}\int_{\partial\Gamma}\phi_m\phi_n\,ds.
\end{equation}
for $1\leq m,n\leq N_{h}$, where $[\cdot]_{mn}$ denotes the $(m,n)$th entry of the corresponding matrix. The vectors $\mathbf p,\mathbf v\in\mathbb R^{N_{h}}$ contain the coefficients of the finite element approximations $p_h$ and $v_h$, respectively.

Hyperbolic wave equations are susceptible to numerical dispersion and artificial dissipation, which may accumulate over long time integrations and introduce numerical artifacts. In uncertainty quantification, such artifacts can be difficult to distinguish from the uncertainty inherent in the inverse problem. To minimize spurious numerical dissipation, we employ the second-order symplectic St\"ormer--Verlet method \cite{hairer2006geometric}, which exhibits excellent long-time stability and accurately preserves the underlying wave dynamics, apart from the physical energy loss induced by the absorbing boundary condition. Applying the St\"ormer--Verlet scheme to the semi-discrete system \eqref{eq:semi-discrete} yields the fully discrete time-stepping scheme
\begin{equation} \label{eq:fully-discrete}
\begin{aligned}
\mathbf p^{n+1/2}
&=
\mathbf p^n+\frac{\Delta t}{2}\mathbf v^n,\\
\mathbf M\mathbf v^{n+1}
&=
\mathbf M\mathbf v^n
-\Delta t\,\mathbf K\mathbf p^{n+1/2}
-\Delta t\,\mathbf B\mathbf v^n,\\
\mathbf p^{n+1}
&=
\mathbf p^{n+1/2}
+\frac{\Delta t}{2}\mathbf v^{n+1}.
\end{aligned}
\end{equation}
The updates for the pressure are explicit, whereas the velocity update requires the application of $\mathbf M^{-1}$ at every time step. With the consistent mass matrix, this amounts to solving a sparse linear system, which dominates the computational cost of the time integration. To obtain a fully explicit and computationally efficient time-stepping scheme, we replace the consistent mass matrix $\mathbf M$ by its lumped approximation \cite{fried1975finite}
\[
\mathbf M \approx \mathbf M_{\mathrm L}
:=
\operatorname{diag}(\mathbf M\mathbf 1),
\]
where $\mathbf 1$ denotes the vector of ones. The resulting diagonal mass matrix can be inverted by componentwise division, eliminating a linear solve at each time step. Mass lumping introduces additional spatial discretization error but provides an efficient explicit approximation of the finite-element wave equation. In our experiments, we observed no qualitative difference between the lumped and consistent solutions at the mesh resolutions considered. More importantly, \Cref{sec:numerical-diagnostics} indicates no discernible bias in posterior uncertainty estimates, while mass lumping reduces the computational cost of the forward and adjoint solves by approximately one order of magnitude. We therefore adopt mass lumping throughout.

Introducing the canonical variable $\mathbf z^n=((\mathbf p^n)^T,(\mathbf v^n)^T)^T$, the St\"ormer--Verlet update \eqref{eq:fully-discrete} can be written as the linear recurrence
\begin{equation} \label{eq:propagator}
    \mathbf z^{n+1}
    =
    \mathbf A_{\Delta t}\mathbf z^n,
\end{equation}
where $\mathbf A_{\Delta t}\in\mathbb R^{2N_{h}\times 2N_{h}}$ is the one-step propagator
\begin{equation} \label{eq:matrix-stromer}
\mathbf A_{\Delta t}
=
\begin{pmatrix}
\mathbf I-\dfrac{\Delta t^2}{2}\mathbf M_{\mathrm L}^{-1}\mathbf K
&
\Delta t\,\mathbf I
-\dfrac{\Delta t^2}{2}\mathbf M_{\mathrm L}^{-1}\mathbf B
-\dfrac{\Delta t^3}{4}\mathbf M_{\mathrm L}^{-1}\mathbf K
\\[2ex]
-\Delta t\,\mathbf M_{\mathrm L}^{-1}\mathbf K
&
\mathbf I
-\Delta t\,\mathbf M_{\mathrm L}^{-1}\mathbf B
-\dfrac{\Delta t^2}{2}\mathbf M_{\mathrm L}^{-1}\mathbf K
\end{pmatrix}.
\end{equation}
Here, $\mathbf I$ denotes the identity matrix of size $N_{h}\times N_{h}$. Although $\mathbf A_{\Delta t}$ is useful for deriving the discrete adjoint, it is never assembled explicitly in the implementation; the action of $\mathbf A_{\Delta t}$ is applied through the three St\"ormer--Verlet substeps in \eqref{eq:fully-discrete}. 

We now derive the discrete counterpart of the forward operator. We assume that ultrasound sensors are located at finite element nodes and, therefore, observation operator reduces to selecting the pressure degrees of freedom. Let $\Pi:\mathbb R^{2N_{h}}\rightarrow\mathbb R^{N_{\mathrm{obs}}}$ denote the projection operator that extracts the pressure values at the sensor locations. Owing to the interpolation property of the Lagrange basis functions,
\begin{equation}\label{eq:lagrange}
    p(\mathbf x_{\mathrm s}^j)
    \approx
    \sum_{i=1}^{N_{\mathrm{FEM}}}
    [\mathbf p]_i\phi_i(\mathbf x_{\mathrm s}^j)
    =
    [\mathbf p]_{k},
\end{equation}
where $k$ is the index of the mesh node coinciding with the $j$th sensor. Hence,
\[
\Pi\mathbf z^n
=
\left[
[\mathbf p^n]_{\mathbb I_1},
\dots,
[\mathbf p^n]_{\mathbb I_{N_{\mathrm{obs}}}}
\right]^T,
\]
where $\mathbb I=\{\mathbb I_1,\dots,\mathbb I_{N_{\mathrm{obs}}}\}$ denotes the sensor node indices. Equivalently, $\Pi$ consists of the corresponding rows of the identity matrix. Like the propagator $\mathbf A_{\Delta t}$, the projection operator is used only for analysis and never assembled explicitly.

Using the discrete propagator, the complete forward operator is given by
\begin{equation}\label{eq:discrete-wave-matrix}
    \mathbf G
    =
    \begin{pmatrix}
        \Pi \mathbf A_{\Delta t}\\
        \Pi \mathbf A_{\Delta t}^{2}\\
        \vdots\\
        \Pi \mathbf A_{\Delta t}^{N_T}
    \end{pmatrix},
\end{equation}
where the exponent denotes matrix powers. The discrete measurements now are

\[
\mathbf y
=
\mathbf G\mathbf z_0,
\]
where $\mathbf z_0=((\mathbf p_0)^T,\mathbf 0^T)^T$, with $\mathbf p_0$ containing the finite element coefficients of the initial pressure distribution and $\mathbf 0\in\mathbb R^{N_h}$ denoting the zero vector. Reshaping the vector $\mathbf y$ into matrix $\mathbf Y$ yields the discrete observation model
\begin{equation}\label{eq:forward-discrete}
    \mathbf Y_{\mathrm{obs}}
    =
    \mathbf Y
    +
    \eta,
\end{equation}
where the additive noise $\eta$ is defined as in \Cref{sec:PA}.

Although the ordering in \eqref{eq:forward-discrete} is convenient for forward simulation, the discrete adjoint derivation benefits from grouping temporal measurements by sensor index. Let $\mathbf R^\pi$ denote this permutation matrix and define the permuted forward operator
\[
\mathbf G^\pi := \mathbf R^\pi \mathbf G.
\]
The corresponding observation model is
\begin{equation}\label{eq:forward-discrete-permuted}
    \begin{pmatrix}
        \mathbf y^1_{\mathrm{obs}}\\
        \vdots\\
        \mathbf y^{N_{\mathrm{obs}}}_{\mathrm{obs}}
    \end{pmatrix}
    =
    \mathbf G^\pi\mathbf z_0
    +
    \eta^\pi,
\end{equation}
where $\eta^\pi=\mathbf R^\pi\eta$ is the permuted noise vector, with unchanged statistical properties. This reordering is used only to derive the adjointeee and is never formed explicitly.
\section{Bayesian Formulation of the PA imaging Problem} \label{sec:Bayes}
In this section, we formulate the PAT inverse problem in the Bayesian framework. The initial pressure and measurements are modeled as random variables related through the discrete forward model. The solution is characterized by the conditional distribution of the initial pressure given the observed data, known as the posterior distribution. This formulation combines prior information with measurement data while accounting for observational uncertainty. For a comprehensive introduction to Bayesian inverse problems, see \cite{kaipio05}.

\subsection{Prior Distribution} \label{sec:prior}
In this work, we assume a Gaussian prior distribution for the discretized initial pressure field,
\[
\mathbf p_0 \sim \mathcal N(\mathbf 0,\mathbf C),
\]
where $\mathbf C\in\mathbb R^{N_{h}\times N_{h}}$ is a symmetric positive-definite covariance matrix. Throughout this work, we adopt a zero-mean prior, although the proposed methodology extends straightforwardly to non-zero prior means. We consider two classes of Gaussian priors: an independent Gaussian prior and the Whittle--Mat\'ern prior \cite{rasmussen2006,roininen2014whittle}. The construction of these covariance models is described in the following subsections.

\subsubsection{Independent Gaussian Prior} \label{sec:white}
A natural first approach is to impose an independent Gaussian prior on the finite element coefficients, treating them as independent standard normal random variables. However, this prior depends on the chosen basis and changes under mesh refinement, tying the induced uncertainty to the discretization rather than the underlying physical field. Since the unknown in PAT is the continuous initial pressure distribution, we instead define the prior with respect to the physical $L^2(\Gamma)$ norm and derive the corresponding distribution of the finite element coefficients. Although the resulting prior remains discretization dependent, its interpretation is consistent with the underlying function-space norm across mesh resolutions. For the independent Gaussian prior, we define the prior density with respect to the physical $L^2(\Gamma)$ norm of the initial pressure field,
\[
    \pi(p_0)
    \propto
    \exp\left(
        -\frac12\|p_0\|_{L^2(\Gamma)}^2
    \right).
\]
Using the expansion $p_0(\boldsymbol x)=\sum_{i=1}^{N_{h}}[\mathbf p_0]_i\phi_i(\boldsymbol x)$, the $L^2(\Gamma)$ norm is approximated by $\|p_0\|_{L^2(\Gamma)}^2\approx \mathbf p_0^T\mathbf M\mathbf p_0$.
In the implementation, consistent with the mass-lumped time discretization, we replace $\mathbf M$ by the lumped mass matrix $\mathbf M_{\mathrm L}=\operatorname{diag}(\mathbf M\mathbf 1)$. Therefore, the discrete prior density is
\begin{equation}\label{eq:iid-prior}
    \pi(\mathbf p_0)
    \propto
    \exp\left(
        -\frac12
        \mathbf p_0^T\mathbf M_{\mathrm L}\mathbf p_0
    \right).
\end{equation}
Equivalently, the finite element coefficient vector satisfies $\mathbf p_0\sim\mathcal N(\mathbf 0,\mathbf M_{\mathrm L}^{-1})$. Since $\mathbf M_{\mathrm L}$ is diagonal, the nodal coefficients are independent Gaussian random variables with variances given by the reciprocal lumped masses,
\[
    [\mathbf p_0]_i
    \sim
    \mathcal N\left(0,\frac{1}{[\mathbf M_{\mathrm L}]_{ii}}\right),
    \qquad i=1,\ldots,N_{\mathrm{FEM}}.
\]
Thus, independence is imposed with respect to the physical $L^2(\Gamma)$ inner product rather than by assigning unit variance directly to the finite element coefficients.

Despite its widespread use in Bayesian PAT, particularly for uncertainty quantification, we do not regard this prior as well suited to large-scale problems. Its coefficient variances increase under mesh refinement, scaling inversely with the local mesh size in the present setting, and may therefore become unrealistically large on fine meshes. Nevertheless, given its prevalent use in Bayesian PAT, we include this prior as a reference for comparison with the alternative prior proposed below.

\subsubsection{Whittle-Mat\'ern Covariance} \label{sec:matern}
The second prior considered in this work is the Whittle--Mat\'ern prior, which introduces spatial correlations through the covariance operator
\begin{equation}\label{eq:matern-cov-op}
    \mathcal C
    =
    \sigma^2 \mathcal Z_{\ell,s}^{-1}
    \left(
        \frac{1}{\ell^2}\mathcal I-\Delta
    \right)^{-(\nu+d/2)},
\end{equation}
where $\sigma$ is the marginal standard deviation, $Z_{\ell,\nu}$ is the normalization constant, $\ell>0$ denotes the correlation length, $d$ is the spatial dimension, and $s$ controls the regularity of the random field. The corresponding prior density is
\begin{equation}\label{eq:prior_model}
    \pi(\mathbf p_0)
    \propto
    \exp\!\left(
        -\frac12
        \mathbf p_0^T
        \mathbf C^{-1}
        \mathbf p_0
    \right),
\end{equation}
where $\mathbf C$ is the finite element approximation of the covariance operator $\mathcal C$.

Typically, $\mathbf C$ is obtained by a finite element discretization of $\mathcal C$. Here, instead, we employ a matrix-free, mesh-independent construction based on a Fourier--spectral representation of the Whittle--Mat'ern field on a box containing the computational domain, followed by interpolation onto the finite element mesh. The correlation parameters $\ell$ and $s$ determine the effective spatial resolution of the field, so the number of dominant modes remains approximately stable under mesh refinement and is typically much smaller than $N_h$.

\subsection{The Likelihood and The Posterior Distributions}
According to the discrete observation model \eqref{eq:forward-discrete-permuted}, the observed data satisfy
\[
\mathbf y_{\mathrm{obs}}^\pi
=
\mathbf G^\pi\mathbf z_0+\eta^\pi,
\qquad
\mathbf z_0=
\begin{pmatrix}
\mathbf p_0\\
\mathbf 0
\end{pmatrix},
\]
where $\mathbf z_0$ is the initial state vector consisting of the finite element coefficients of the initial pressure and the zero initial velocity. We assume a Gaussian likelihood with known absolute noise level $\sigma:=d_{\mathrm{noise}}\frac{\|\mathbf Y\|_{L^2(0,T;\ell^2)}}{\|\eta_0\|_{L^2(0,T;\ell^2)}}$, according to \eqref{eq:noise-model}.
%,where $d_{\mathrm{noise}}$ is the prescribed relative noise level and $\eta_0$ is the standard Gaussian realization used to generate the synthetic measurements. 
The likelihood is therefore given by
\begin{equation}
\pi(\mathbf y_{\mathrm{obs}}^\pi\mid\mathbf p_0)
\propto
\exp\left(
-\frac{1}{2\sigma^2}
\|
\mathbf G^\pi\mathbf z_0-\mathbf y_{\mathrm{obs}}^\pi
\|_{L^2(0,T;\ell^2)}^2
\right).
\end{equation}
Since $\mathbf y^\pi$ is obtained by a permutation of the entries of $\mathbf Y$, the discrete $L^2(0,T;\ell^2)$ norm is preserved under the permutation and is therefore used interchangeably for both representations. Combining the likelihood with the Gaussian prior introduced in \Cref{sec:prior} using Bayes' theorem gives the posterior density
\begin{equation}\label{eq:posterior}
\pi(\mathbf p_0\mid\mathbf y_{\mathrm{obs}}^\pi)
\propto
\exp\left(
-\frac{1}{2\sigma^2}
\|
\mathbf G^\pi\mathbf z_0-\mathbf y_{\mathrm{obs}}^\pi
\|_{L^2(0,T;\ell^2)}^2
-\frac12\mathbf p_0^T\mathbf C^{-1}\mathbf p_0
\right).
\end{equation}
We remark that the synthetic noise realization is normalized to attain the prescribed relative noise level $d_{\mathrm{noise}}$ exactly. This defines the absolute noise scale $\sigma$ used in the likelihood. Consequently, the synthetic perturbation is not strictly Gaussian. In the Bayesian formulation, however, $\sigma$ is assumed to be known and fixed, and we adopt the corresponding Gaussian likelihood model. Since the forward operator is linear and the prior is Gaussian, the resulting posterior distribution is also Gaussian.

\section{The Randomized-Then-Optimize Method For Exploring The Posterior} \label{sec:RTO}
The maximum a posteriori (MAP) estimate is obtained by maximizing the posterior density \eqref{eq:posterior}, or equivalently, by minimizing its negative logarithm. Since the forward model is linear and both the likelihood and the prior are Gaussian, the MAP estimate is the unique solution of the quadratic optimization problem
\begin{equation}\label{eq:MAP}
    \mathbf p_0^{\mathrm{MAP}}
    =
    \underset{\mathbf p_0}{\operatorname{arg\,min}}
    \;
    \frac{1}{\sigma^2}
    \left\|
        \mathbf G^\pi
        \begin{pmatrix}
            \mathbf p_0\\
            \mathbf 0
        \end{pmatrix}
        -
        \mathbf y_{\mathrm{obs}}^\pi
    \right\|_{L^2(0,T;\ell^2)}^2
    +
    \mathbf p_0^T\mathbf C^{-1}\mathbf p_0.
\end{equation}
The randomized-then-optimized (RTO) method \cite{bardsley2014randomize} generates posterior samples by solving randomly perturbed versions of the MAP problem. For the linear-Gaussian inverse problem considered here, each posterior sample is obtained as
\begin{equation} \label{eq:optimization}
    \mathbf p_0^{\mathrm{sample}}
    =
    \underset{\mathbf p_0}{\operatorname{arg\,min}}
    \;
    \frac{1}{\sigma^2}
    \left\|
        \mathbf G^\pi
        \begin{pmatrix}
            \mathbf p_0\\
            \mathbf 0
        \end{pmatrix}
        -
        \left(
            \mathbf y_{\mathrm{obs}}^\pi
            +
            \sigma \boldsymbol\eta^\pi
        \right)
    \right\|_{L^2(0,T;\ell^2)}^2
    +
    \left\|
        \mathbf C^{-1/2}
        \left(
            \mathbf p_0
            -
            \mathbf p_{\mathrm{prior}}
        \right)
    \right\|_2^2,
\end{equation}
where $[\boldsymbol\eta^\pi]_{j,i}\stackrel{\mathrm{i.i.d.}}{\sim}\mathcal N(0,1)$ and $\mathbf p_{\mathrm{prior}}=\mathbf C^{1/2}\boldsymbol\xi$ with $\boldsymbol\xi\sim\mathcal N(\mathbf 0,\mathbf I)$.

The RTO method transforms the deterministic MAP problem into a stochastic optimization problem through the random perturbations in \eqref{eq:optimization}. For linear inverse problems with Gaussian likelihood and prior, each perturbed solution is a posterior sample \cite{bardsley2014randomize,wang2017bayesian}. Thus, independent perturbations yield independent posterior samples without constructing a Markov chain, reducing posterior sampling to randomized least-squares problems.

The RTO formulation presented in \eqref{eq:optimization} differs slightly from the original formulation proposed by \cite{bardsley2014randomize}. The original method is based on a stacked least-squares formulation, which is mathematically elegant but less naturally suited to large-scale matrix-free finite element implementations. In contrast, the formulation adopted here preserves the original forward and adjoint operators and introduces randomness solely through the right-hand side of the optimization problem derived below. Consequently, the same matrix-free operator can be reused for every posterior sample, with only the right-hand side changing between optimization problems. To demonstrate the equivalence of the two formulations in the linear-Gaussian setting, we revisit the RTO construction of \cite{bardsley2014randomize} in Appendix~\ref{appendix:rto-equivalence}.

To formally express the solution to the randomized optimization problem \eqref{eq:optimization}, we first rewrite it as an equivalent least-squares problem. Defining the extended operator and the right-hand side
\begin{equation} \label{eq:extended-matrices}
    \mathbf G_{\mathrm{ext}}
    :=
    \begin{pmatrix}
        \sigma^{-1}\mathbf G^\pi\\
        \mathbf M_{\mathrm L}^{1/2}
    \end{pmatrix}, \qquad 
    \mathbf b
    :=
    \begin{pmatrix}
        \sigma^{-1}
        \left(
            \mathbf y_{\mathrm{obs}}^\pi
            +
            \sigma\boldsymbol\eta^\pi
        \right)\\
        \mathbf M_{\mathrm L}^{1/2}
        \mathbf p_{\mathrm{prior}}
    \end{pmatrix},
\end{equation}
the optimization problem \eqref{eq:optimization} can be written as
\begin{equation}
    \mathbf p_0^{\mathrm{sample}}
    =
    \underset{\mathbf p_0}{\operatorname{arg\,min}}
    \;
    \|
    \mathbf G_{\mathrm{ext}}\mathbf p_0
    -
    \mathbf b
    \|_2^2.
\end{equation}
The corresponding normal equations are
\begin{equation}
\label{eq:normal-equations}
    \mathbf G_{\mathrm{ext}}^T
    \mathbf G_{\mathrm{ext}}
    \mathbf p_0^{\mathrm{sample}}
    =
    \mathbf G_{\mathrm{ext}}^T
    \mathbf b.
\end{equation}
Iterative least-squares methods, such as LSQR \cite{paige1982lsqr} or CGLS \cite{hansen1998rank}, solve \eqref{eq:normal-equations} without explicitly forming the normal matrix or the extended operator $\mathbf G_{\mathrm{ext}}$. Here, neither $\mathbf G^\pi$ nor $\mathbf G_{\mathrm{ext}}$ is assembled. Instead, the action of $\mathbf G^\pi$ is evaluated matrix-free by repeatedly applying the one-step propagator $\mathbf A_{\Delta t}$ and observation operator $\Pi$, as described in Section~\ref{sec:FEM}. The transpose $\mathbf G_{\mathrm{ext}}^T$, required by LSQR, is derived next.

\subsection{Matrix-free transpose and the continuous adjoint} \label{sec:implicit-transpose}
The transpose operator $\mathbf G_{\mathrm{ext}}^T$ plays a central role in the matrix-free solution of the normal equations \eqref{eq:normal-equations}. While the action of the forward operator $\mathbf G$ is readily evaluated by repeated application of the one-step propagator introduced in Section~\ref{sec:FEM}, the corresponding transpose operator is not immediately available. Since neither $\mathbf G$ nor $\mathbf G_{\mathrm{ext}}$ is assembled explicitly, the transpose must also be realized in a matrix-free manner.

Constructing such a transpose operator presents both conceptual and computational challenges. From a computational perspective, the transpose cannot be obtained by explicitly transposing the forward matrix, as this matrix is never formed. From a mathematical perspective, one may either derive the transpose of the fully discrete forward operator or discretize the continuous adjoint wave equation. These two procedures are not generally equivalent, since discretization and adjoint formation do not commute. This distinction is particularly important for iterative least-squares methods such as LSQR, whose convergence theory relies on the availability of the exact algebraic transpose of the discrete forward operator.

In this section we first characterize the transpose of the fully discrete forward operator by transposing the one-step propagator. We then identify the corresponding continuous evolution equation and compare it with the classical continuous adjoint formulation for the photoacoustic wave equation. This comparison reveals a discrepancy arising from the artificial boundary treatment. Finally, we show that by enlarging the computational domain sufficiently, the influence of this discrepancy can be confined outside the region of interest during the observation interval, thereby reconciling the computational and physical formulations.

Recall that the discrete forward operator is obtained by repeatedly applying the one-step propagator $\mathbf A_{\Delta t}$ in \eqref{eq:propagator} together with the observation operator $\mathbf y^n=\Pi\mathbf z^n$. Let $\widehat{\mathbf y}=((\widehat{\mathbf y}^{\,1})^T,\ldots,(\widehat{\mathbf y}^{\,N_T})^T)^T$ be an arbitrary vector, where $\widehat{\mathbf y}^{\,n}\in\mathbb R^{N_{\mathrm s}}$. By the definition of the forward operator,
\begin{equation} \label{eq:transpose-definition}
    \mathbf G^T\widehat{\mathbf y}
    =
    \sum_{n=1}^{N_T}
    (\mathbf A_{\Delta t}^{\,n})^T
    \Pi^T
    \widehat{\mathbf y}^{\,n}.
\end{equation}
Using the identity $(\mathbf A_{\Delta t}^{\,n})^T=((\mathbf A_{\Delta t})^T)^n$, the action of $\mathbf G^T$ can be evaluated recursively without explicitly forming powers of $\mathbf A_{\Delta t}$. To this end, define $\widehat{\mathbf z}^{\,N_T} = \Pi^T\widehat{\mathbf y}^{\,N_T}$, and compute recursively, for $n=N_T-1,\ldots,1$,
\begin{equation}
\label{eq:backward-recursion}
    \widehat{\mathbf z}^{\,n}
    =
    \Pi^T\widehat{\mathbf y}^{\,n}
    +
    (\mathbf A_{\Delta t})^T
    \widehat{\mathbf z}^{\,n+1}.
\end{equation}
We remark that, since the sensors are assumed to coincide with finite element mesh nodes, the operator $\Pi^T$ simply injects the measurement residuals into the corresponding degrees of freedom. In more general settings, where the sensors do not coincide with mesh nodes, $\Pi^T$ is the transpose of the interpolation operator and may distribute the residuals over neighboring basis functions. Expanding the recursion shows that
\[
    (\mathbf A_{\Delta t})^T\widehat{\mathbf z}^{\,1}
    =
    (\mathbf A_{\Delta t})^T\Pi^T\widehat{\mathbf y}^{\,1}
    +
    (\mathbf A_{\Delta t}^{\,2})^T\Pi^T\widehat{\mathbf y}^{\,2}
    +\cdots+
    (\mathbf A_{\Delta t}^{\,N_T})^T\Pi^T\widehat{\mathbf y}^{\,N_T},
\]
which is precisely $\mathbf G^T\widehat{\mathbf y}$. Therefore, the action of the transpose operator is obtained by a backward recursion involving only the transpose of a time step. The remaining task is to characterize the transpose of the one-step St\"ormer--Verlet propagator $(\mathbf A_{\Delta t})^T$. 

The transpose of the one-step propagator $\mathbf A_{\Delta t}$ is obtained by taking the algebraic transpose of the block matrix in \eqref{eq:propagator} with derivations provided in \Cref{sec:appendix-transppose}. Introducing the adjoint variables $\widehat{\mathbf p}$ and $\widehat{\mathbf v}$, one backward time step consists of the following operations:
\begin{equation}\label{eq:adjoint-verlet}
\begin{aligned}
\widehat{\mathbf p}^{\,n+\frac12}
&=
\widehat{\mathbf p}^{\,n+1},\\
\widehat{\mathbf v}^{\,n+\frac12}
&=
\widehat{\mathbf v}^{\,n+1}
+\frac{\Delta t}{2}\widehat{\mathbf p}^{\,n+1},\\
\mathbf M_{\mathrm L}\widehat{\mathbf q}^{\,n+\frac12}
&=
\widehat{\mathbf v}^{\,n+\frac12},\\
\widehat{\mathbf p}^{\,n}
&=
\widehat{\mathbf p}^{\,n+\frac12}
-\Delta t\,\mathbf K^T\widehat{\mathbf q}^{\,n+\frac12},\\
\widehat{\mathbf v}^{\,n}
&=
\left(\mathbf M_{\mathrm L}-\Delta t\,\mathbf B^T\right)
\widehat{\mathbf q}^{\,n+\frac12}
+\frac{\Delta t}{2}\widehat{\mathbf p}^{\,n}.
\end{aligned}
\end{equation}
In this backward time-marching scheme, the substeps are applied in reverse order, with each linear operator replaced by its algebraic transpose. Consequently, the transpose operator can be evaluated in a matrix-free manner using the same sparse operators employed by the forward solver. In the present work, the stiffness and boundary matrices, $\mathbf K$ and $\mathbf B$, are symmetric, so their transpose coincides with the original operator. More generally, however, the transpose of a sparse matrix is readily available and can be applied with essentially the same computational cost, so the proposed framework extends naturally to non-symmetric discretizations.

The variables appearing in \eqref{eq:adjoint-verlet} are algebraic transpose variables associated with the Euclidean inner product on the finite element coefficient space. In contrast, finite element discretizations of the wave equation are naturally expressed with respect to the $L^2$ inner product. We therefore introduce the mass-dual variables $\boldsymbol\lambda_p^n=\mathbf M_{\mathrm L}^{-1}\widehat{\mathbf p}^{\,n}$ and $\boldsymbol\lambda_v^n=\mathbf M_{\mathrm L}^{-1}\widehat{\mathbf v}^{\,n}$ which may be interpreted as the finite element coefficients of the adjoint pressure and velocity fields. Expressing the backward recursion in these transformed variables applies the discrete Riesz map from the coefficient-space dual variables to finite element coefficient vectors (see \Cref{sec:appendix-andjoint-continuous} for more detail). Consequently, the recursion assumes the standard Galerkin form $\mathbf M_{\mathrm L} \dot{\boldsymbol\lambda} + \mathbf K \boldsymbol{\lambda}$, making it possible to identify the corresponding continuous first-order wave system.

Formally, ignoring the artificial boundary contribution, the corresponding continuous adjoint dynamics are
\begin{equation}\label{eq:continuous-adjoint-formal}
\begin{aligned}
    \frac{\partial \lambda_v}{\partial \tau} - \Delta \lambda_p
    &=
    -\sum_{n=1}^{N_T}
    \sum_{i=1}^{N_{\mathrm s}}
    [\widehat{\mathbf y}^{\,n}]_i
    \delta_{\mathbf x_{\mathrm s}^i}(\mathbf x)
    \delta(\tau-t_n),
    && \mathbf x\in\Gamma,\quad \tau\in[0,T],\\
    \frac{\partial \lambda_p}{\partial \tau}-\lambda_v
    &=0,
    && \mathbf x\in\Gamma,\quad \tau\in[0,T],
\end{aligned}
\end{equation}
Here $\tau=T-t$ denotes reverse time. The reverse-time wave equation \eqref{eq:continuous-adjoint-formal} is the continuous adjoint of the PAT forward model on an unbounded domain. The algebraic transpose of the fully discrete solver differs from this adjoint only through the artificial boundary contribution introduced by domain truncation. We therefore extend the computational domain sufficiently so that any boundary reflections do not re-enter the physical domain during $[0,T]$. This confines boundary artifacts outside the region of interest and motivates our boundary treatment over a perfectly matched layer (PML). In our numerical experiments, the PML implementation introduced artifacts in the posterior uncertainty estimates, further motivating this choice.

To summarize the forward and transpose operators introduced in this section, recall that each LSQR iteration consists of one application of the extended forward operator followed by one application of its transpose. Given an iterate $\mathbf p$, the extended forward action is
\begin{equation} \label{eq:lsqr-residual}
\mathbf r
=
\mathbf G_{\mathrm{ext}}\mathbf p
-
\mathbf b
=
\begin{pmatrix}
\sigma^{-1}
\left(
\mathbf G^\pi
\begin{pmatrix}
\mathbf p\\
\mathbf 0
\end{pmatrix}
-
\mathbf y_{\mathrm{obs}}^\pi
-
\sigma\boldsymbol\eta^\pi
\right)
\\[0.4em]
\mathbf M_{\mathrm L}^{1/2}
\left(
\mathbf p-\mathbf p_{\mathrm{prior}}
\right)
\end{pmatrix},
\end{equation}
where $\mathbf r=((\mathbf r_{\mathrm d})^T,(\mathbf r_{\mathrm p})^T)^T$ denotes the residual of the augmented least-squares problem. The corresponding transpose action is
\begin{equation} \label{eq:lsqr-transpose-action}
\widehat{\mathbf p}
=
\mathbf G_{\mathrm{ext}}^T\mathbf r
=
\sigma^{-1}
\mathbf J^T(\mathbf G^\pi)^T
\mathbf r_{\mathrm d}
+
\mathbf M_{\mathrm L}^{1/2}
\mathbf r_{\mathrm p},
\qquad
\mathbf J\mathbf p=
\begin{pmatrix}
\mathbf p\\
\mathbf 0
\end{pmatrix}.
\end{equation}
The remaining LSQR operations consist of vector updates and scalar normalizations.

In the present implementation, the first block of the forward residual is evaluated by one forward St\"ormer--Verlet simulation together with the observation operator, while the transpose action is computed by the backward recursion \eqref{eq:backward-recursion} and the transpose St\"ormer--Verlet scheme \eqref{eq:adjoint-verlet}. Consequently, each LSQR iteration requires exactly one forward wave propagation, one reverse-time propagation, and a small number of sparse vector operations. No global matrices are assembled or factorized. The procedure for drawing an RTO posterior sample is summarized in \Cref{alg:rto-sample}.

\begin{algorithm}[t]
\caption{One RTO posterior sample}
\label{alg:rto-sample}
\begin{algorithmic}[1]
\State Draw $\boldsymbol\eta^\pi\sim\mathcal N(\mathbf0,\mathbf I)$ and $\mathbf p_{\mathrm{prior}}\sim\mathcal N(\mathbf0,\mathbf C)$.
\State Define $\mathbf b$ as in \eqref{eq:extended-matrices}.
\State Use LSQR to solve $\min_{\mathbf p}\|\mathbf G_{\mathrm{ext}}\mathbf p-\mathbf b\|_2^2$.
\While{LSQR stopping criterion is not satisfied}
    \State Apply the forward St\"ormer--Verlet solver to compute $\mathbf G^\pi(\mathbf p^T,\mathbf0^T)^T$.
    \State Assemble the augmented residual $\mathbf r=(\mathbf r_{\mathrm d}^T,\mathbf r_{\mathrm p}^T)^T$ as in \eqref{eq:lsqr-residual}.
    \State Apply the backward recursion \eqref{eq:backward-recursion} and transpose St\"ormer--Verlet scheme \eqref{eq:adjoint-verlet} to compute $(\mathbf G^\pi)^T\mathbf r_{\mathrm d}$.
    \State Assemble $\widehat{\mathbf p}=\mathbf G_{\mathrm{ext}}^T\mathbf r$ as in \eqref{eq:lsqr-transpose-action} using the previous step.
    \State Perform the LSQR vector updates and normalizations.
\EndWhile
\State \Return $\mathbf p^{\mathrm{sample}}=\mathbf p_{\mathrm{LSQR}}$.
\end{algorithmic}
\end{algorithm}

\section{Numerical Results} \label{sec:results}
In this section, we present numerical experiments that validate the proposed Bayesian formulation and demonstrate its performance for two- and three-dimensional PAT problems. We first report numerical verification tests, including an explicit posterior variance computation on a small mesh, LSQR convergence behavior, and the effect of mass lumping on accuracy and computational cost. We then consider two two-dimensional test cases. The first is a disk phantom, where we compare reconstructions and posterior uncertainty for different sensor coverage lengths, including full, half, and quarter boundary measurements. The second is a limited-angle blood-vessel phantom motivated by experimental data, for which we compare the proposed RTO sampler with NUTS. Finally, we present three-dimensional examples with both independent Gaussian and Whittle--Mat\'ern priors to demonstrate the applicability of the method to large-scale 3D PAT reconstruction. The accompanying source code is publicly available on GitHub. \cite{afkham2026pat3duq}.

\subsection{Numerical Diagnostics} \label{sec:numerical-diagnostics}
Before presenting the reconstruction results, we perform a series of numerical diagnostics to validate the proposed methodology. These include verification of the matrix-free transpose operator, a reference posterior covariance on a small-scale problem, an investigation of the convergence behavior of LSQR for both 2D and 3D problems, and an assessment of the accuracy and computational efficiency of the mass-lumped discretization.

To verify the proposed transpose operator, we performed the discrete adjoint test $\langle \mathbf G\mathbf p,\widehat{\mathbf y}\rangle_{\ell^2}=\langle \mathbf p,\mathbf G^T\widehat{\mathbf y}\rangle_{\ell^2}$ using independently generated vectors $\mathbf p$ and $\widehat{\mathbf y}$ with i.i.d. standard normal entries. In all experiments, the relative discrepancy $\bigl|\langle \mathbf G\mathbf p,\widehat{\mathbf y}\rangle_{\ell^2}-\langle \mathbf p,\mathbf G^T\widehat{\mathbf y}\rangle_{\ell^2}\bigr|/\bigl|\langle \mathbf G\mathbf p,\widehat{\mathbf y}\rangle_{\ell^2}\bigr|$ was on the order of machine precision, confirming that the backward recursion \eqref{eq:backward-recursion} together with \eqref{eq:adjoint-verlet} implements the exact algebraic transpose of the discrete forward operator.

Since the proposed Bayesian formulation differs from standard PAT formulations through the space--time noise model introduced in \Cref{eq:noise-model}, we first construct a reference posterior solution. To this end, we explicitly assemble the discrete forward operator $\mathbf G$ on a small mesh by exploiting the identity $\mathbf G\mathbf e_i=[\mathbf G]_i$, where $\mathbf e_i$ is the $i$th Euclidean basis vector and $[\mathbf G]_i$ denotes the $i$th column of $\mathbf G$. The assembled matrix allows the posterior covariance $\mathbf C_{\mathrm{post}}=\left(\frac{1}{\sigma^2}\mathbf G^T\mathbf G+\mathbf C^{-1}\right)^{-1}$ to be computed explicitly, from which the pointwise posterior variance is obtained as the diagonal of the covariance matrix. Throughout this experiment, we employ a two-dimensional mesh with $20\times20$ interior grid points and a time step of $\Delta t\approx 0.028$. The phantom used in this model is given by
\begin{equation} \label{eq:disk-phantom-formula}
p_0(x,y)
=
{1}/(1+\exp\!\left(\left(\sqrt{(x-0.5)^2+(y-0.5)^2}-0.18\right)/0.015\right)).
\end{equation}

\begin{figure}[t]
    \centering
    \includegraphics[width=\linewidth]{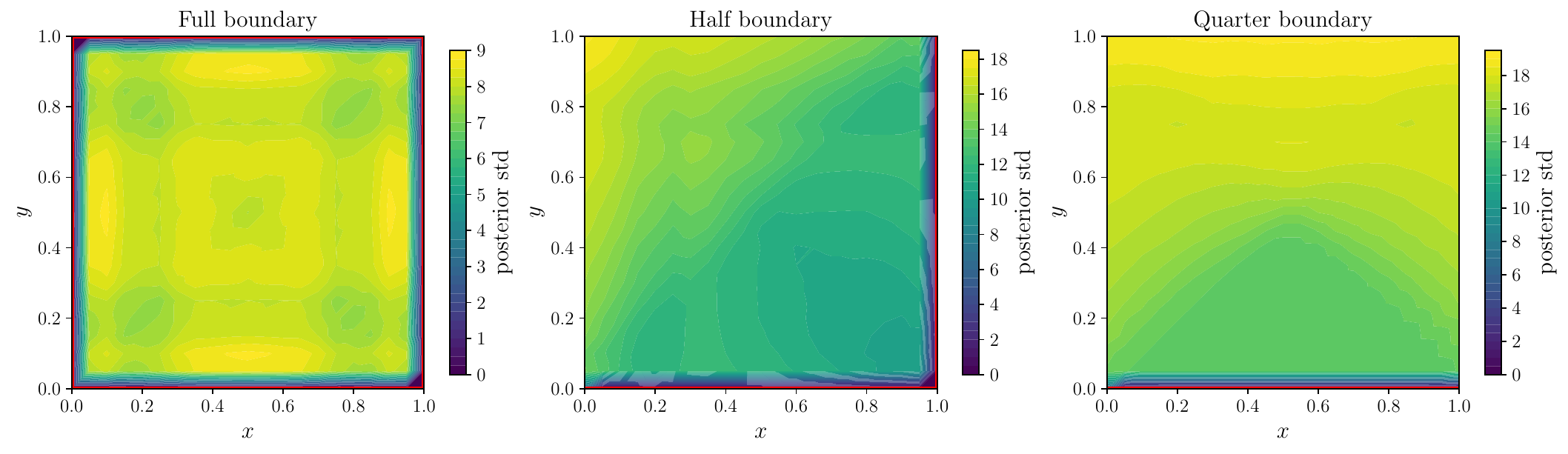}
    \caption{
    Reference pointwise posterior variance computed from the explicitly assembled posterior covariance on a $20\times20$ interior mesh. Red segments indicate the measurement boundary.
    }
    \label{fig:ref-solution}
\end{figure}

We present the reference pointwise posterior variance in \Cref{fig:ref-solution}. As expected, the posterior uncertainty decreases as the measurement aperture increases. For partial boundary measurements, the uncertainty grows with increasing distance from the measurement boundary. This behavior is consistent with the Bayesian uncertainty estimates reported in \cite{tarvainen2013bayesian}. We emphasize that the explicit computation of the posterior covariance is only feasible for relatively small problems, since it requires assembling the discrete forward operator and forming the dense posterior covariance matrix, whose computational cost and memory requirements become prohibitive as the number of unknowns increases.

\begin{figure}[t]
    \centering
    \includegraphics[width=0.6\linewidth]{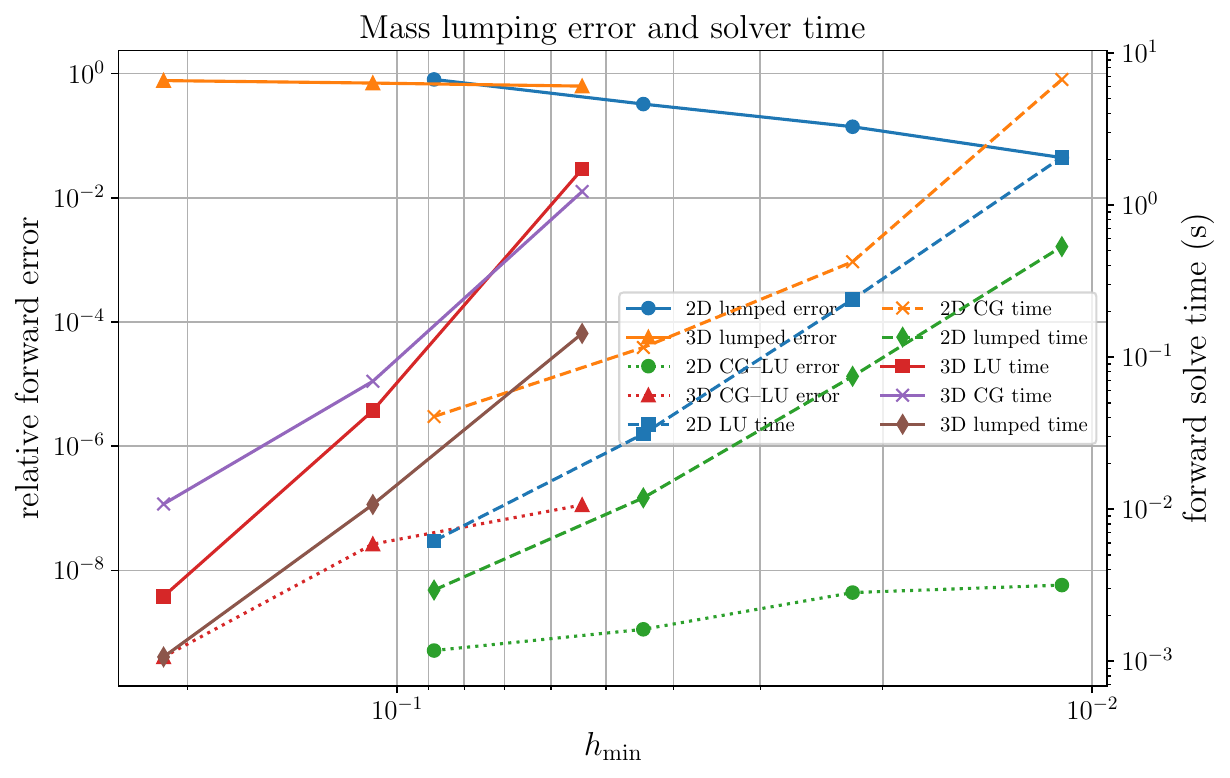}
    \caption{
    Relative forward error (left vertical axis) and computation time (right vertical axis) for the mass-lumped, conjugate gradient (CG), and LU-factorization approaches for solving the mass matrix system $\textbf M \textbf x=\textbf b$.
    }
    \label{fig:mass-lumping}
\end{figure}

We present the mass-lumping error for the 2D and 3D test cases in \Cref{fig:mass-lumping}, where the mass-lumped approximation is compared with LU factorization and the conjugate gradient (CG) method for solving $\mathbf M\mathbf x=\mathbf b$. Although the approximation error is non-negligible, it decreases under mesh refinement and can therefore be controlled by increasing the mesh resolution, with a more pronounced effect observed in two dimensions. In return, mass lumping provides up to an order-of-magnitude reduction in computational time compared with LU factorization and approximately a fivefold speedup over CG. Furthermore, we did not observe any systematic differences in the posterior uncertainty estimates, and therefore employ mass lumping throughout the remaining numerical experiments.

\begin{figure}[t]
    \centering
    \includegraphics[width=\linewidth]{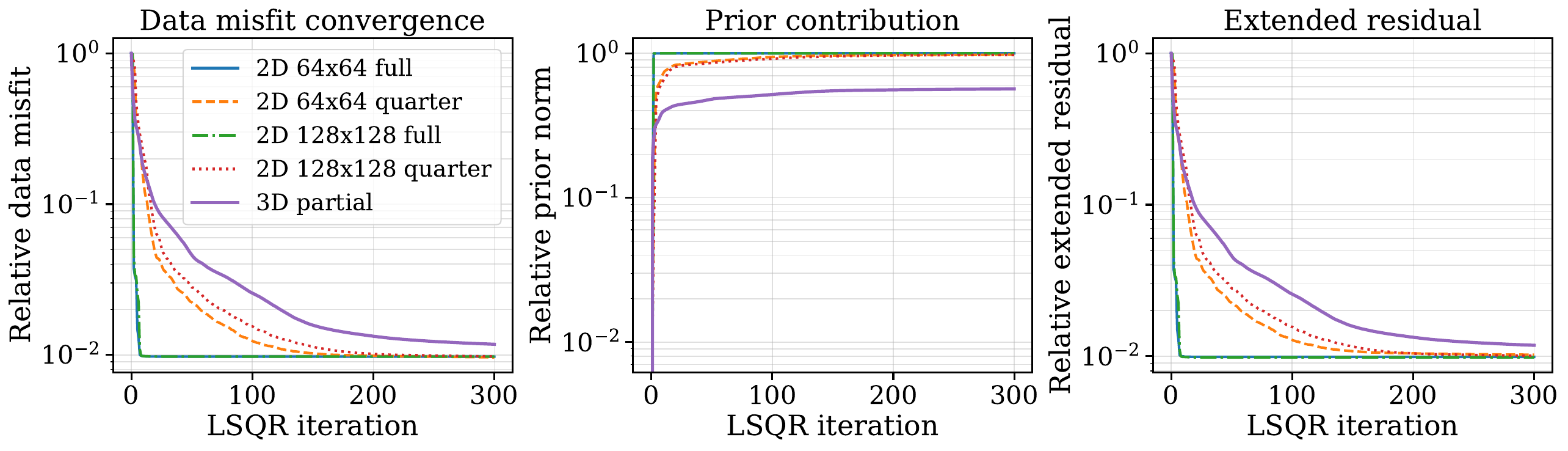}
    \caption{
    LSQR convergence diagnostics for the two- and three-dimensional test cases. The normalized data misfit, prior contribution, and extended residual are shown as functions of the iteration count.
    }
    \label{fig:lsqr-diagnostics}
\end{figure}

Finally, we present representative LSQR convergence diagnostics for the 2D problems with full and partial boundary measurements on $128\times128$ and $64\times64$ FEM meshes, and for the 3D problem with partial boundary measurements in \Cref{fig:lsqr-diagnostics}. We report the normalized data residual, normalized prior residual, and normalized residual of the augmented least-squares system. A commonly used stopping criterion for iterative regularization methods is the discrepancy principle, the whitened data residual satisfies $\|\mathbf G\mathbf p-\mathbf y\|_{L^2(0,T;\ell^2)}\approx\sigma\sqrt{N_T\times N_{\mathrm{obs}}}$,
which is typically achieved within only $2$--$10$ LSQR iterations for all problems considered. While this number of iterations is sufficient to obtain an accurate MAP estimate and posterior mean, it is insufficient for reliable uncertainty quantification, as significantly more iterations are required to obtain representative pointwise posterior standard deviation estimates. On the other hand, owing to the well-known semi-convergence behavior of Krylov subspace methods, excessive iterations eventually introduce noise-dominated components that deteriorate both the reconstruction and the estimated posterior uncertainty. Ideally, the parameter space should be decomposed into a data-informed subspace, which is solved using LSQR, and a complementary prior-dominated subspace, which is sampled directly from the prior without further iterations. Developing such a decomposition and the corresponding stopping strategy is beyond the scope of the present work and is left for future research. Consequently, throughout this work we employ fixed iteration limits of $300$ LSQR iterations for the 2D experiments and $350$ LSQR iterations for the 3D experiments. The choice of stopping criterion for RTO in severely ill-posed PDE-constrained inverse problems remains an important practical and theoretical question.

\subsection{2D PAT Problem}
In this section, we present the numerical results for the 2D PAT problem. We first consider a disk phantom with full, half, and quarter boundary measurements to compare the estimated posterior uncertainty with the reference covariance presented in the previous section and with previously reported results in the literature. We then investigate a limited-view blood vessel phantom and compare the proposed RTO sampler with the No-U-Turn Sampler (NUTS), a variant of Hamiltonian Monte Carlo.

\begin{figure}[t]
    \centering
    \includegraphics[width=\linewidth]{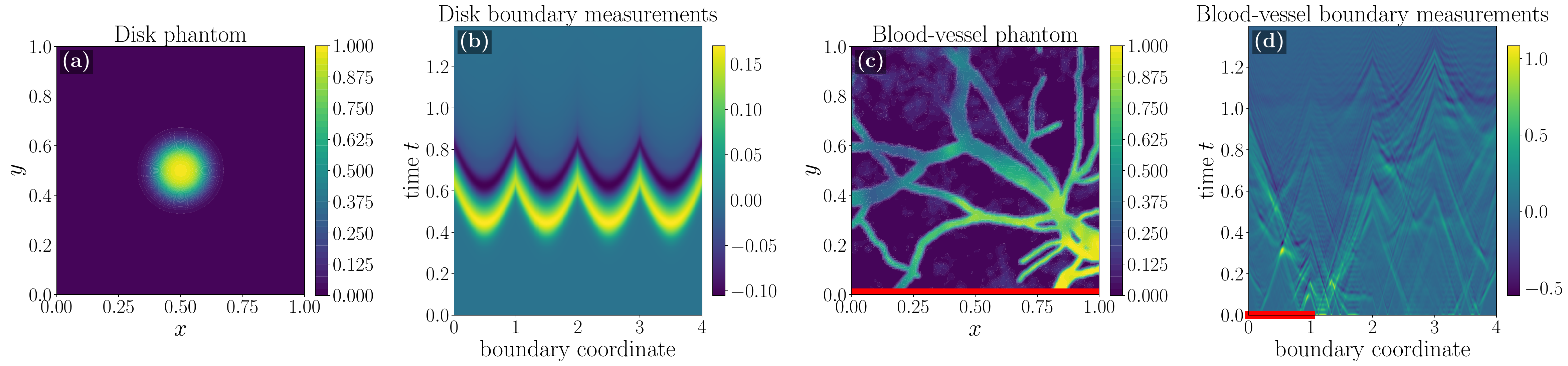}
    \caption{
    True initial pressures and corresponding boundary measurements for the 2D disk and blood-vessel phantoms. The measured boundary segment is highlighted in red.
    }
    \label{fig:ground-truth-2D}
\end{figure}

\subsubsection{Disk Phantom}

All two-dimensional experiments in this section are performed on the unit square $\Gamma=[0,1]^2$, discretized using a structured triangular mesh with $64\times64$ interior grid points. To mitigate the influence of the artificial boundary while avoiding the additional memory required to keep the boundary inactive throughout the entire experiment, the computational domain is enlarged by $0.75$ units in every direction using the same mesh spacing. Since acoustic information leaves the physical domain after approximately $T=1.4$, this extension ensures that any waves reflected from the artificial boundary cannot return to the physical domain during the observation interval. The wave equation is discretized using the mass-lumped St\"ormer--Verlet scheme introduced in \Cref{sec:FEM} with $\Delta t=0.0012$, satisfying the CFL stability condition. Throughout this section, we employ the independent Gaussian prior introduced in \Cref{sec:prior}, corresponding to the precision matrix $\mathbf M_{\mathrm L}$.

The initial pressure is generated from \eqref{eq:disk-phantom-formula} and is provided in \Cref{fig:ground-truth-2D}. To create the measurements we solve the forward problem on a mesh of size $128\times128$ which coincides with the computational mesh and time steps. Pressure variations are then collected at every mesh node on the computational mesh on the active portion of the boundary and at every time step. This ensures that we avoid interpolation error between the two meshes. We set the relative noise level to $d_{\mathrm{noise}}=0.01$ and repeat the experiment for three choices of $\Gamma_{\mathrm{obs}}$: the full boundary, one half of the boundary, and one quarter of the boundary.

For each case, we construct the posterior distribution as described in \Cref{sec:Bayes} and generate posterior samples using the RTO procedure in \Cref{alg:rto-sample}. Given samples $\{\mathbf p_0^{(k)}\}_{k=1}^{N_{\mathrm{sample}}}$, the pointwise posterior mean and standard deviation are estimated by
\begin{equation}\label{eq:ergodic}
\overline{\mathbf p}_0
=
\frac{1}{N_{\mathrm{sample}}}
\sum_{k=1}^{N_{\mathrm{sample}}}
\mathbf p_0^{(k)},
\qquad
\boldsymbol\sigma_{\mathbf p_0}
=
\left[
\frac{1}{N_{\mathrm{sample}}-1}
\sum_{k=1}^{N_{\mathrm{sample}}}
\left(
\mathbf p_0^{(k)}-\overline{\mathbf p}_0
\right)^2
\right]^{1/2},
\end{equation}
where the square and square root are taken component-wise. For each measurement configuration, we generate $N_{\text{sample}}=100$ posterior samples. Although this sample size is modest, it provides a stable estimate of the pointwise posterior variance owing to the fact that RTO produces independent posterior samples. In contrast, Markov chain Monte Carlo methods such as NUTS generate correlated samples, typically requiring substantially larger sample sizes to overcome autocorrelation and achieve a comparable effective sample size.

\begin{figure}[t]
    \centering
    \includegraphics[width=\linewidth]{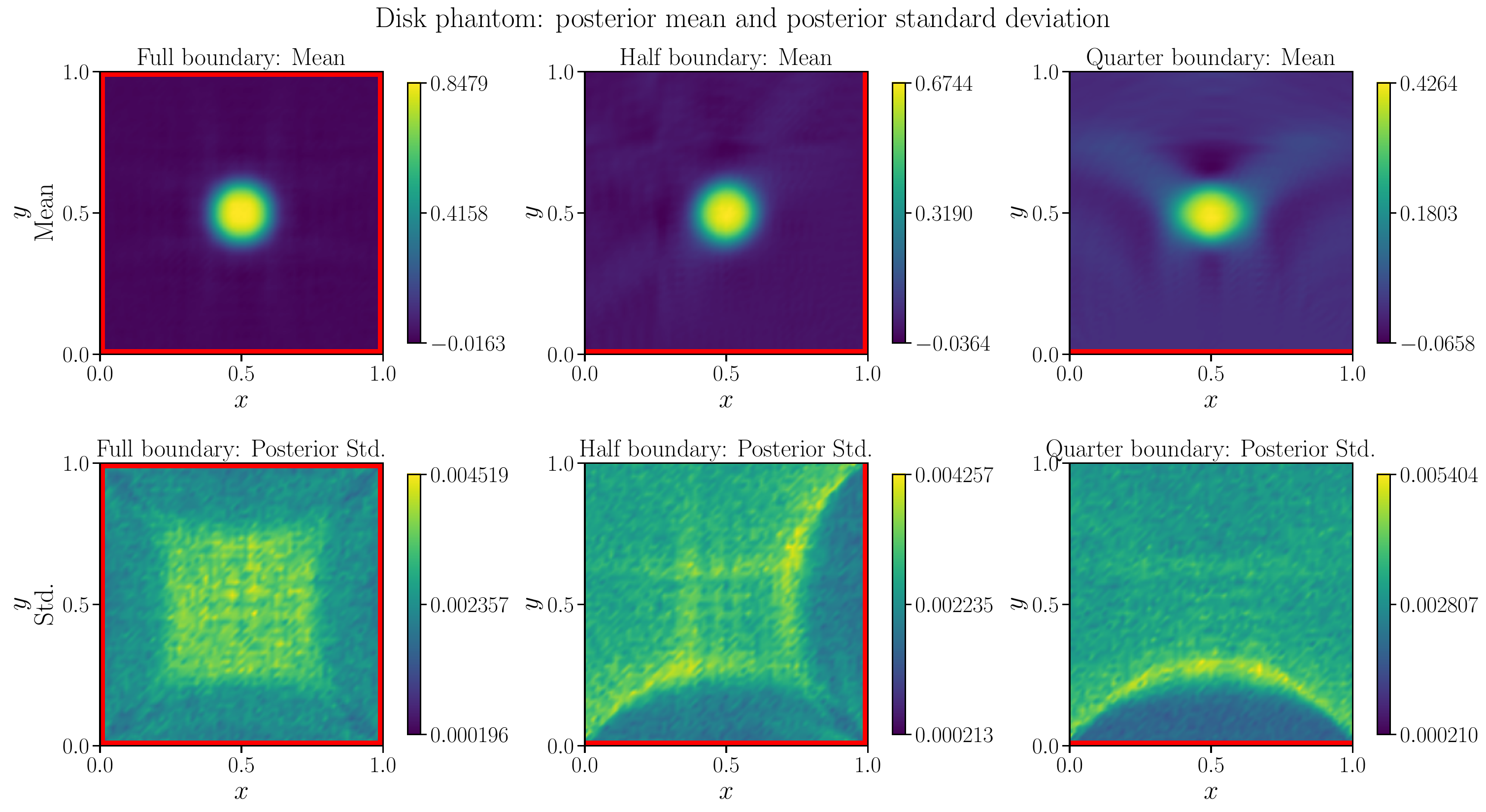}
    \caption{
    Empirical mean and pointwise posterior standard deviation estimated from RTO posterior samples for the disk phantom. Red segments indicate the measurement boundary. The red lines indicate the measurement boundaries.
    }
    \label{fig:disk-results}
\end{figure}

We present the mean estimate, which coincides with the posterior MAP for the linear Gaussian inverse problem considered here, together with the pointwise posterior standard deviation in \Cref{fig:disk-results}. The mean reconstructions obtained from the full- and half-boundary measurements accurately recover the location and shape of the inclusion, whereas the quarter-boundary configuration exhibits the expected degradation associated with limited-view data. These reconstruction characteristics are consistent with previously reported results for limited-view PAT \cite{tarvainen2013bayesian,tick2016}. The estimated posterior standard deviations also agree qualitatively with the reference covariance in \Cref{fig:ref-solution}, up to a normalization constant, exhibiting increasing uncertainty with distance from the measurement boundary. Compared with the reference solution, the sample-based standard deviation estimates exhibit a moderate level of Monte Carlo noise, which can be reduced by increasing the number of posterior samples. Since the numerical approximations introduced throughout this work can be systematically controlled through mesh refinement, solver tolerances, and additional sampling, the quality of the estimated posterior statistics can be improved accordingly.

Overall, the sample-based posterior statistics are in good qualitative agreement with the explicitly computed reference covariance while remaining computationally feasible on meshes for which explicit covariance construction is no longer practical. In addition, the posterior standard deviation is smallest near the measurement boundary and increases with distance from the measurement aperture, consistent with the behavior reported for limited-view PAT.

\subsection{Blood Vessels Phantom in a Limited View}

In this section, we evaluate the proposed method on the blood vessels phantom shown in \Cref{fig:ground-truth-2D}(c) and compare the resulting posterior statistics with those obtained using the NUTS Sampler, a variant of Hamiltonian Monte Carlo. We restrict our attention to the quarter-boundary measurement configuration (indicated by the red line), as this represents the most ill-posed setting considered in this work and therefore provides the most informative assessment of the uncertainty estimates.

The computational setup follows that of the disk phantom experiment. The physical domain is discretized using a structured $128\times128$ mesh and extended by $0.75$ units in each coordinate direction using the same mesh spacing to form the computational domain. This yields a finite element discretization with 103041 degrees of freedom. The time step is chosen as $\Delta t=0.001$, satisfying the CFL stability condition.

The spatio-temporal measurement data are presented in \Cref{fig:ground-truth-2D}(d). In the partial-view configuration, the measurements are restricted to the portion of the boundary indicated by the red line. Throughout this experiment, the relative noise level is set to $d_{\mathrm{noise}} = 0.01$.

We compare the performance of RTO with the No-U-Turn Sampler (NUTS), an adaptive variant of Hamiltonian Monte Carlo (HMC). NUTS augments the posterior distribution with an auxiliary momentum variable and employs Hamiltonian dynamics to efficiently explore the posterior distribution. As a gradient-based Markov chain Monte Carlo method, each proposal requires evaluation of the gradient of the negative log-posterior, which in PAT is obtained by solving the corresponding adjoint problem described above. We select NUTS as a benchmark because it is among the most successful and widely used gradient-based sampling algorithms for Bayesian inverse problems. We used Pyro, a Python package, for all NUTS implementations. In this example we consider 100 RTO samples and 1200 NUTS samples.

\begin{figure}[t]
    \centering
    \includegraphics[width=\linewidth]{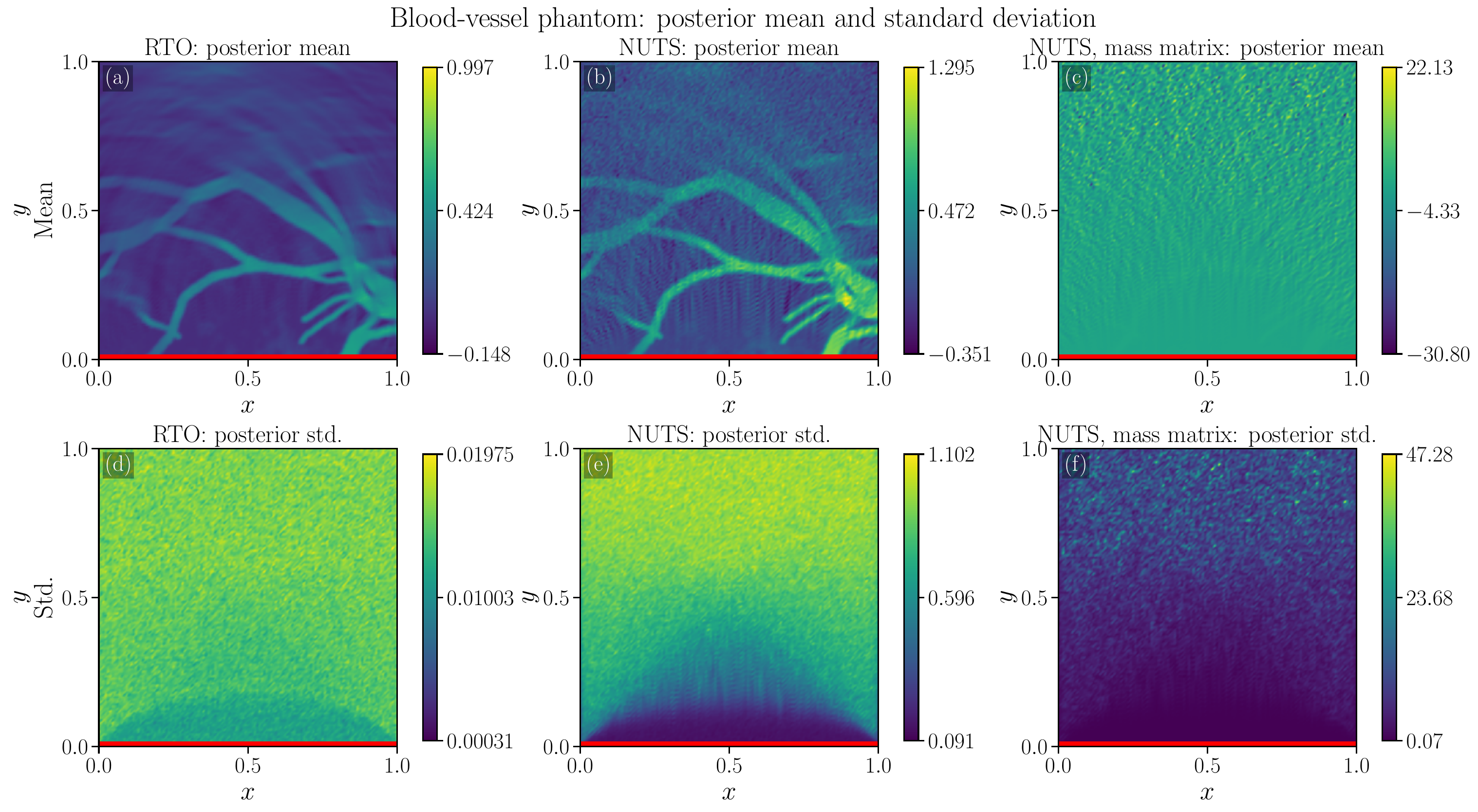}
    \caption{
    Empirical mean and pointwise posterior standard deviation estimated from RTO and NUTS posterior samples for the blood vessels phantom. The red lines indicate the measurement boundaries.
    }
    \label{fig:blood-vessels-results}
\end{figure}

We present the posterior mean and pointwise standard deviation obtained with RTO in \Cref{fig:blood-vessels-results}(a) and (d), respectively. As in the disk phantom, the uncertainty exhibits the expected spatial behavior. For comparison, the corresponding NUTS estimates are shown in \Cref{fig:blood-vessels-results}(e) and (f). Despite the Gaussian posterior, NUTS fails to adequately explore the distribution. A possible explanation is that standard NUTS samples the initial momentum in the FEM coefficient space, where the mass-lumping prior can map many momentum realizations to nearly zero physical configurations. To investigate this effect, we also consider a standard Gaussian prior, $\mathbf p_0 \sim \mathcal N(0,I_{N_h})$, obtained by removing the mass-lumping prior. The corresponding posterior mean and standard deviation are shown in \Cref{fig:blood-vessels-results}(b) and (e). While this modification improves the sampling performance of NUTS, it degrades the posterior mean because of the nonphysical prior. Moreover, the estimated standard deviation and mean exhibits artifacts, including diagonal and vertical streaks, that are likely induced by the nonphysical prior.

We note a discrepancy between the magnitudes of the posterior standard deviations estimated by NUTS and RTO, consistent with the observations for the disk phantom. A possible explanation is the LSQR stopping criterion, which may terminate the iterations either too early or too late.

The sample quality was poor for both NUTS implementations. We ran four chains, each with 300 warm-up iterations followed by 300 posterior samples, yielding a total of 1200 posterior samples. For the standard Gaussian prior, the effective sample size was only 91, which is substantially lower than the performance typically achieved by NUTS on differentiable inverse problems. This result highlights the difficulty of sampling high-dimensional PDE-constrained posterior distributions.

In summary, RTO provides an embarrassingly parallel approach for posterior sampling that naturally accommodates FEM and other PDE-based constraints arising in PAT while accurately capturing the spatial structure of the posterior standard deviation. However, a stopping criterion for LSQR that accounts for both the prior and the likelihood is needed to accurately recover its magnitude. These experiments also highlight limitations of NUTS. Although recent advances have extended NUTS to infinite-dimensional settings, further work is needed to make these methods effective for PDE-constrained inverse problems with FEM discretizations.

\subsection{3D PAT Problem}

\begin{figure}[t]
    \centering
    \includegraphics[width=\linewidth]{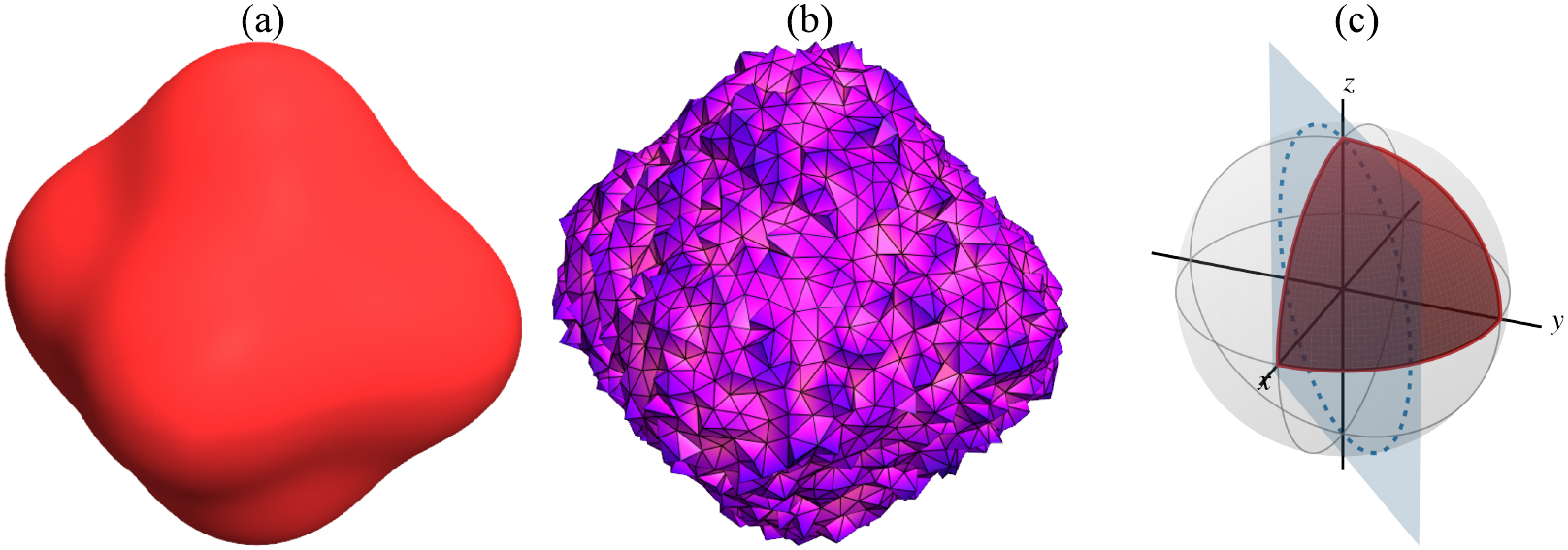}
    \caption{
    True phantom (only the inclusion surface) and computational geometry for the 3D PAT problem. (a) ground truth inclusion. (b) Projection of the inclusion onto the finite element mesh. (c) Physical domain showing the measurement boundary (red) and the visualization plane. 
    }
    \label{fig:ground-truth-3D}
\end{figure}

In this section, we investigate a three-dimensional limited-view PAT problem to demonstrate the applicability of the proposed Bayesian framework to large-scale geometries. The physical domain is chosen to be the unit sphere, while the measurement boundary is restricted to the portion of the surface lying in the first octant, $\Gamma_{\mathrm{obs}}=\{\mathbf{x}\in\partial\Omega:\;x>0,\;y>0,\;z>0\}$. The true initial pressure is represented by a smooth non-spherical inclusion obtained by perturbing a sphere with a low-order spherical-harmonic-type angular deformation and positioned slightly off the center of the domain. It is defined by
\begin{equation}
p_0(\mathbf{x})
=
\left(
1+\exp\left(
\frac{
\|\mathbf{x}-\mathbf{c}\|_2
-
R_0\!\left(1+\alpha S\!\left(\dfrac{\mathbf{x}-\mathbf{c}}{\|\mathbf{x}-\mathbf{c}\|_2}\right)\right)
}{\varepsilon}
\right)
\right)^{-1}.
\end{equation}
Here, $R_0=0.5$, $\alpha=0.8$, $\varepsilon=0.03$, and
$\mathbf{c}=\left(0.1,0.1,0.1\right)/\sqrt{3}$.
The function $S$ denotes a normalized quartic angular perturbation that produces a smooth non-spherical inclusion. Consequently, $p_0$ represents an approximately constant background initial pressure with an elevated pressure inside the inclusion. The true phantom, its finite element approximation, and the computational domain together with the measurement geometry are shown in Fig.~\ref{fig:ground-truth-3D}.

The computational domain is obtained by extending the physical domain by one unit in every direction, resulting in a spherical computational domain of radius $2$. As in the two-dimensional experiments, this extension ensures that acoustic waves leave the physical domain before any reflections from the artificial boundary can return during the observation interval. The computational domain is discretized using a uniform unstructured tetrahedral mesh with 207151 degrees of freedom. The mesh is generated such that the measurement surface is represented exactly, ensuring that all measurement locations coincide with mesh nodes and eliminating interpolation errors. The wave equation is discretized using the mass-lumped St\"ormer-Verlet scheme introduced in \Cref{sec:FEM} with a time step of $\Delta t\approx0.0097$, satisfying the CFL stability condition. The final simulation time is chosen as $T=2$, which is sufficient for all acoustic waves to leave the physical domain. The measurement data are generated by solving the forward problem on the same computational mesh, and a relative noise level of $d_{\mathrm{noise}}=0.01$ is employed throughout this section.

In this section, we consider two prior distributions. The first is the independent Gaussian prior introduced in \Cref{sec:prior}, $\mathbf p_0 \sim \mathcal N(\mathbf 0,\mathbf M_{\mathrm L}^{-1})$. The second is a Whittle-Mat\'ern prior with correlation length $\ell=0.15$, smoothness parameter $s=1.5$, marginal standard deviation $\sigma=3.0$, and zero mean. To generate samples from this prior, we employ the covariance operator \eqref{eq:matern-cov-op} and solve the stochastic partial differential equation $\left(\frac{1}{0.15^2}\mathcal I-\Delta\right)^3p_0=W$, where $W$ denotes Gaussian white noise. Here we have used $d=3$, so that the exponent satisfies $s+d/2=3$. The equation is solved using a Fourier spectral method on the cube $[-2.5,2.5]^3$ with a $128^3$ uniform grid. The resulting samples are then projected onto the finite element space defined on the spherical computational mesh.

For both prior models, posterior samples are generated using the RTO algorithm described in \Cref{alg:rto-sample}. As in the two-dimensional experiments, we generate $100$ independent posterior samples by solving the corresponding augmented least-squares problems using LSQR with a fixed iteration limit of $350$. In the independent Gaussian case, the unknown is represented directly in the finite element basis, whereas for the Whittle--Mat\'ern prior the RTO optimization is carried out in the latent Gaussian parameter space associated with the Fourier representation of the prior. Each latent sample is subsequently mapped to the finite element space before evaluating the forward model.

Since the RTO algorithm is embarrassingly parallel, we employed a parallel implementation using 10 message passing interface (MPI) processes, with each process generating 10 posterior samples. The complete computation required approximately 5 hours. Because the posterior samples are independent, the workload scales almost linearly with the number of available processes. Consequently, if each process were assigned a single sample, the total computation time would be reduced to approximately 30. All experiments were conducted on a MacBook Pro equipped with an M4 Max processor and 48 GB of memory.

\begin{figure}[t]
    \centering
    \includegraphics[width=\linewidth]{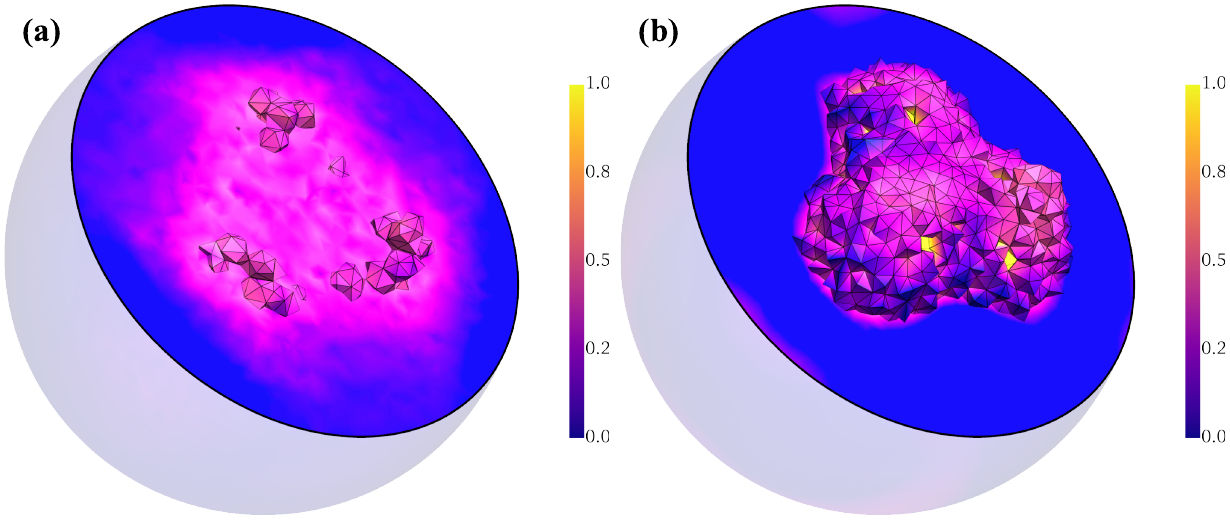}
    \caption{
    Mean/MAP estimates for the posterior constructed with (a) the independent Gaussian prior and (b) the Whittle-Matérn prior.  
    }
    \label{fig:posterior-mean-3D}
\end{figure}

We illustrate the posterior mean, computed using the ergodic averages \eqref{eq:ergodic}, for both the independent Gaussian prior \eqref{eq:ergodic}(a) and the Whittle--Matérn prior \eqref{eq:ergodic}(b) in \Cref{fig:posterior-mean-3D}. Each panel displays a planar cross-section of the physical domain together with the three-dimensional inclusion obtained by thresholding the posterior mean at the value (0.5). The cross-sectional plane has normal vector proportional to $(1,1,1)$ and is located at a distance of approximately $0.35$ from the origin.

The independent Gaussian prior recovers the overall geometry of the inclusion but underestimates the initial pressure within the inclusion and produces a more diffuse transition between the inclusion and the background. In contrast, the Whittle-Matérn prior yields a substantially more accurate reconstruction, with a sharper boundary separating the inclusion from the background and a more accurate estimate of the initial pressure inside the inclusion.ion.

\begin{figure}[t]
    \centering
    \includegraphics[width=\linewidth]{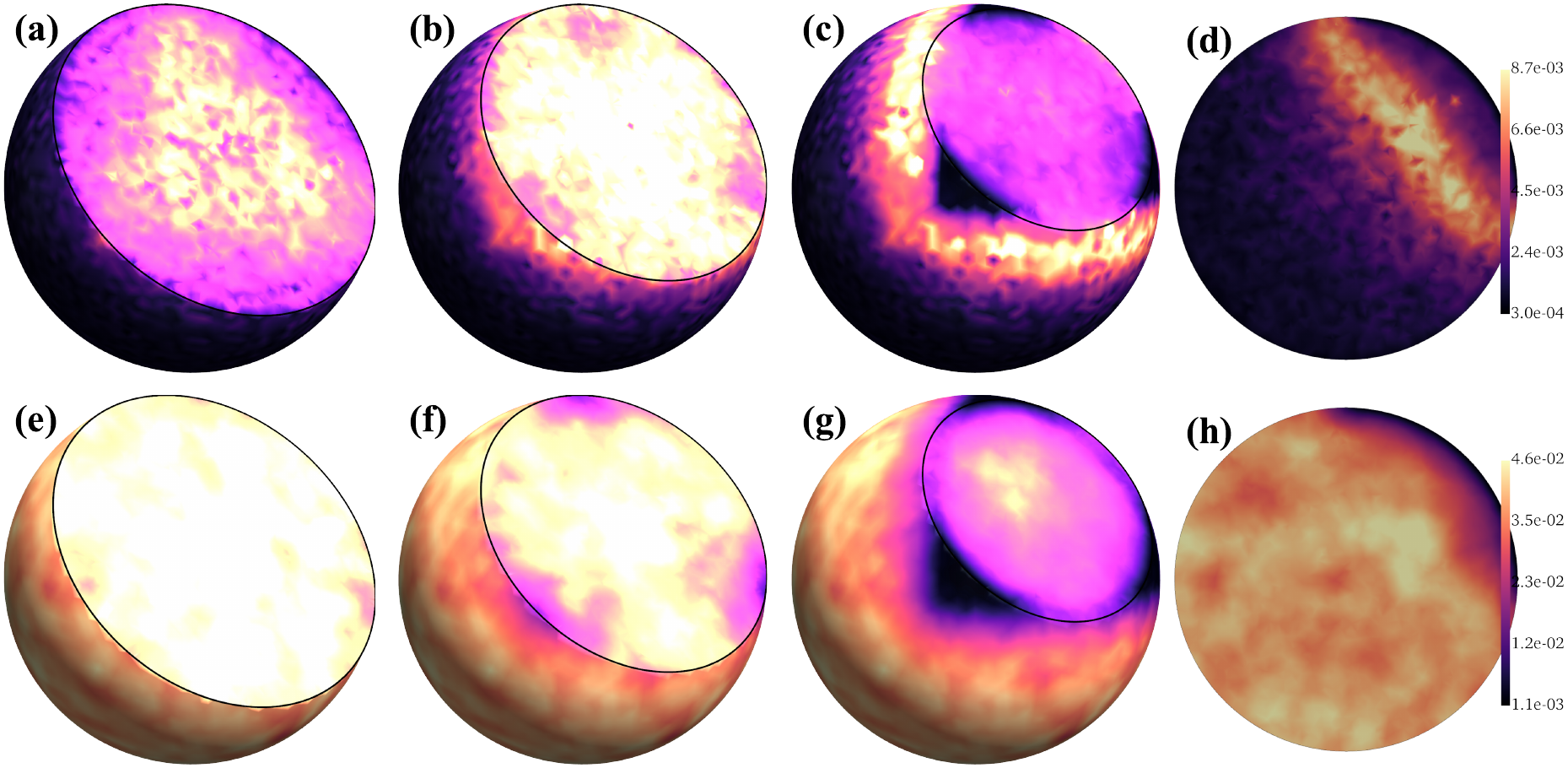}
    \caption{
    Point-wise standard deviation. The first and second rows correspond to the independent Gaussian prior and the Whittle-Matérn prior, respectively. The first three columns show cross-sections of the computational domain at different distances from the centre, while the last column displays the cross-section along the plane illustrated in \Cref{fig:ground-truth-3D}.
    }
    \label{fig:posterior-std-3D}
\end{figure}

We present the poin-twise posterior standard deviation for the independent Gaussian prior (first row) and the Whittle-Matérn prior (second row) in \Cref{fig:posterior-std-3D}. The first three columns show planar slices of the posterior standard deviation obtained using planes with normal vector proportional to $(1,1,1)$ located at signed distances $0$, $0.35$, and $0.75$ from the origin, respectively. These planes move progressively closer to the measurement surface. The fourth column shows the planar cross-section illustrated in \Cref{fig:ground-truth-3D}, whose normal is proportional to $(1,-1,0)$ and which bisects the measurement surface into two equal parts.

For the independent Gaussian prior, the posterior uncertainty decreases as the cross-sectional plane approaches the measurement surface, reflecting the increased information content of the data in this region. However, regions of the domain furthest from the measurement surface also exhibit relatively low posterior standard deviation, which is contrary to the expected behavior for a limited-view inverse problem. A plausible explanation is the early termination of the LSQR iterations. In large-scale and severely ill-posed problems, the initial Krylov subspaces generated by LSQR are dominated by directions corresponding to the largest singular values of the forward operator and are therefore primarily informed by the data. Directions that are weakly constrained by the observations typically enter the Krylov subspace only in later iterations. Since the independent Gaussian prior has a flat covariance spectrum and is not dimension independent, there is no dominant low-dimensional prior subspace, so a substantially larger number of Krylov iterations may be required before these weakly informed directions are adequately represented.

In contrast, the Whittle-Matérn prior is dimension independent and possesses a covariance operator whose eigenvalues decay with increasing frequency. Consequently, the prior information is concentrated in a comparatively small number of dominant modes, allowing the Krylov iterations to capture the most influential prior directions much more rapidly. As shown in the second row of \Cref{fig:posterior-std-3D}, the posterior uncertainty decreases near the measurement surface while remaining appreciably larger in regions that are poorly observed. This behavior is more consistent with the limited-view geometry of the inverse problem and suggests that the Whittle-Matérn prior provides a more faithful quantification of the reconstruction uncertainty.

\section{Conclusions} \label{sec:conclusion}
In this work, we developed a finite-element Bayesian uncertainty quantification framework for photoacoustic tomography that naturally accommodates complex computational domains and detector geometries. The proposed formulation combines a matrix-free randomize-then-optimize (RTO) sampling strategy with a discrete adjoint construction that forms an exact transpose pair with the forward operator while remaining consistent with the continuous PAT adjoint. This enables efficient large-scale posterior sampling using LSQR without assembling the forward operator or explicitly constructing the posterior covariance. Through a series of numerical diagnostics and two- and three-dimensional examples, we demonstrated that the proposed methodology accurately reproduces reference posterior statistics where these are available and extends Bayesian uncertainty quantification to realistic large-scale finite-element PAT problems.

The numerical experiments further demonstrate that reliable Bayesian uncertainty quantification depends not only on the sampling algorithm but also on the numerical treatment of the underlying wave equation. In particular, temporal discretization, artificial boundary conditions, adjoint consistency, iterative solver convergence, and prior modeling all influence the computed posterior uncertainty and should therefore be regarded as integral components of the uncertainty quantification procedure. Although the proposed methodology provides an effective framework for large-scale Bayesian PAT, an important open problem is the development of principled stopping criteria for the LSQR iterations within the RTO algorithm. Unlike deterministic inverse problems, where stopping criteria are typically designed to recover the MAP estimate, posterior sampling additionally requires accurate representation of posterior uncertainty, leading to a different trade-off between computational cost and statistical accuracy. Finally, although developed in the context of PAT, the proposed framework extends naturally to a broad class of large-scale linear PDE-constrained inverse problems.

\appendix

\section{Equivalence of the Matrix-Free RTO Formulation}
\label{appendix:rto-equivalence}

In this appendix, we reformulate the RTO method of \cite{bardsley2014randomize} for large-scale matrix-free finite element computations. The equivalent formulation preserves the original forward and adjoint operators, allowing the same matrix-free implementation for all posterior samples while avoiding explicit construction and orthogonalization of the stacked least-squares operator.

Let $\mathbf p\in\mathbb R^n$ denote the unknown parameter and let
$\mathbf y\in\mathbb R^m$ denote the observed data. We consider the
linear observation model $\mathbf y = \mathbf G\mathbf p + \boldsymbol\eta$, with $\boldsymbol\eta\sim\mathcal N(\mathbf 0,\sigma^2\mathbf I)$, together with the Gaussian prior $\mathbf p\sim\mathcal N(\mathbf 0,\mathbf C)$. The negative log posterior, up to an additive constant, is then
\begin{equation}\label{eq:appendix-negative-log-posterior}
    \Phi(\mathbf p)
    =
    \frac{1}{2\sigma^2}
    \|\mathbf G\mathbf p-\mathbf y\|_2^2
    +
    \frac12
    \mathbf p^T\mathbf C^{-1}\mathbf p .
\end{equation}
Equivalently, if $\mathbf L$ is such that $\mathbf L^T\mathbf L=\mathbf C^{-1}$, then
\[
    \Phi(\mathbf p)
    =
    \frac{1}{2\sigma^2}
    \|\mathbf G\mathbf p-\mathbf y\|_2^2
    +
    \frac12
    \|\mathbf L\mathbf p\|_2^2 .
\]
Since the forward operator $\mathbf G$ is linear and both the likelihood and prior are Gaussian, the posterior distribution is also Gaussian. In particular, $\mathbf p \mid \mathbf y \sim \mathcal N(\boldsymbol\mu_{\mathrm{post}},\mathbf C_{\mathrm{post}})$, where the posterior covariance and mean are given by
\begin{equation}
    \mathbf C_{\mathrm{post}}
    =
    \left(
        \frac{1}{\sigma^2}\mathbf G^T\mathbf G
        +
        \mathbf C^{-1}
    \right)^{-1},
    \qquad
    \boldsymbol\mu_{\mathrm{post}}
    =
    \mathbf C_{\mathrm{post}}
    \frac{1}{\sigma^2}
    \mathbf G^T\mathbf y.
\end{equation}
We now consider the randomly perturbed optimization problem introduced in
Section~\ref{sec:RTO}. For each posterior sample, we draw independent
perturbations
\[
    \boldsymbol\eta \sim \mathcal N(\mathbf 0,\mathbf I),
    \qquad
    \mathbf p_{\mathrm{prior}}\sim\mathcal N(\mathbf 0,\mathbf C),
\]
and solve
\begin{equation}\label{eq:appendix-rto-minimization}
    \mathbf p^{\mathrm{sample}}
    =
    \underset{\mathbf p}{\operatorname{arg\,min}}
    \;
    \frac{1}{2\sigma^2}
    \left\|
        \mathbf G\mathbf p
        -
        \left(\mathbf y+\sigma\boldsymbol\eta\right)
    \right\|_2^2
    +
    \frac12
    \left\|
        \mathbf L
        \left(\mathbf p-\mathbf p_{\mathrm{prior}}\right)
    \right\|_2^2 .
\end{equation}
The minimizer of \eqref{eq:appendix-rto-minimization} is obtained from the corresponding normal equations,
\[
    \left(
        \frac{1}{\sigma^2}\mathbf G^T\mathbf G
        +
        \mathbf L^T\mathbf L
    \right)
    \mathbf p^{\mathrm{sample}}
    =
    \frac{1}{\sigma^2}\mathbf G^T
    \left(
        \mathbf y+\sigma\boldsymbol\eta
    \right)
    +
    \mathbf L^T\mathbf L\mathbf p_{\mathrm{prior}} .
\]
Where the inverse posterior covariance appears on the left-hand-side. Thus,
\begin{equation}\label{eq:appendix-rto-solution}
    \mathbf p^{\mathrm{sample}}
    =
    \mathbf C_{\mathrm{post}}
    \left[
        \frac{1}{\sigma^2}\mathbf G^T\mathbf y
        +
        \frac{1}{\sigma}\mathbf G^T\boldsymbol\eta
        +
        \mathbf L^T\mathbf L\mathbf p_{\mathrm{prior}}
    \right].
\end{equation}
It can be seen that $\mathbb E[\mathbf p^{\mathrm{sample}}] = \boldsymbol\mu_{\mathrm{post}}$,since $\mathbb E[\boldsymbol\eta]=\mathbf 0$ and $\mathbb E[\mathbf p_{\mathrm{prior}}]=\mathbf 0$. Furthermore,
\begin{align}
    \operatorname{Cov}(\mathbf p^{\mathrm{sample}})
    &=
    \mathbf C_{\mathrm{post}}
    \operatorname{Cov}\left(
        \frac{1}{\sigma}\mathbf G^T\boldsymbol\eta
        +
        \mathbf L^T\mathbf L\mathbf p_{\mathrm{prior}}
    \right)
    \mathbf C_{\mathrm{post}} \nonumber\\
    &=
    \mathbf C_{\mathrm{post}}
    \left(
        \frac{1}{\sigma^2}\mathbf G^T\mathbf G
        +
        \mathbf L^T\mathbf L
    \right)
    \mathbf C_{\mathrm{post}} \nonumber=
    \mathbf C_{\mathrm{post}}
\end{align}
where we have used the independence of $\boldsymbol\eta$ and
$\mathbf p_{\mathrm{prior}}$, together with
$\operatorname{Cov}(\boldsymbol\eta)=\mathbf I$ and $ \operatorname{Cov}(\mathbf p_{\mathrm{prior}}) = \mathbf C =(\mathbf L^T\mathbf L)^{-1}$.
Therefore, $ \mathbf p^{\mathrm{sample}} \sim \mathcal N(\boldsymbol\mu_{\mathrm{post}},\mathbf C_{\mathrm{post}})$,
showing that the solution of the perturbed optimization problem is an exact posterior sample.

\section{Transpose St\"ormer--Verlet Propagator} \label{sec:appendix-transppose}
The one-step St\"ormer--Verlet propagator can be written as
\[
\mathbf A_{\Delta t}
=
\mathbf S_4\mathbf S_3\mathbf S_2\mathbf S_1,
\]
where
\[
\mathbf S_1=
\begin{pmatrix}
\mathbf I & \dfrac{\Delta t}{2}\mathbf I\\
\mathbf0&\mathbf I
\end{pmatrix},
\qquad
\mathbf S_2=
\begin{pmatrix}
\mathbf I&\mathbf0\\
-\Delta t\,\mathbf K&
\mathbf M_{\mathrm L}-\Delta t\,\mathbf B
\end{pmatrix},
\]
\[
\mathbf S_3=
\begin{pmatrix}
\mathbf I&\mathbf0\\
\mathbf0&\mathbf M_{\mathrm L}^{-1}
\end{pmatrix},
\qquad
\mathbf S_4=
\begin{pmatrix}
\mathbf I & \dfrac{\Delta t}{2}\mathbf I\\
\mathbf0&\mathbf I
\end{pmatrix}.
\]

Consequently,
\[
\mathbf A_{\Delta t}^T
=
\mathbf S_1^T
\mathbf S_2^T
\mathbf S_3^T
\mathbf S_4^T,
\]
where
\[
\mathbf S_1^T=
\begin{pmatrix}
\mathbf I&\mathbf0\\
\dfrac{\Delta t}{2}\mathbf I&\mathbf I
\end{pmatrix},
\qquad
\mathbf S_2^T=
\begin{pmatrix}
\mathbf I&
-\Delta t\,\mathbf K^T\\
\mathbf0&
(\mathbf M_{\mathrm L}-\Delta t\,\mathbf B)^T
\end{pmatrix},
\]
\[
\mathbf S_3^T=
\begin{pmatrix}
\mathbf I&\mathbf0\\
\mathbf0&\mathbf M_{\mathrm L}^{-T}
\end{pmatrix},
\qquad
\mathbf S_4^T=
\begin{pmatrix}
\mathbf I&\mathbf0\\
\dfrac{\Delta t}{2}\mathbf I&\mathbf I
\end{pmatrix}.
\]
Therefore, the transpose propagator is obtained by applying the transpose of each elementary substep in reverse order. Writing these four matrix-vector products as update equations yields the backward time-marching scheme \eqref{eq:adjoint-verlet}.

\section{Interpretation of the Transpose Variables} \label{sec:appendix-andjoint-continuous}
The variables appearing in the transpose recursion \eqref{eq:adjoint-verlet} are algebraic transpose variables associated with the Euclidean inner product on the finite element coefficient space. To recover the finite element coefficients of the physical adjoint fields, we apply the discrete $L^2$ Riesz map
$\boldsymbol{\lambda}_p^n
=
\mathbf M_{\mathrm L}^{-1}
\widehat{\mathbf p}^{\,n},
$ and $
\boldsymbol{\lambda}_v^n
=
\mathbf M_{\mathrm L}^{-1}
\widehat{\mathbf v}^{\,n}.
$
Equivalently,
$
\widehat{\mathbf p}^{\,n}
=
\mathbf M_{\mathrm L}\boldsymbol{\lambda}_p^n,
$ and $
\widehat{\mathbf v}^{\,n}
=
\mathbf M_{\mathrm L}\boldsymbol{\lambda}_v^n.
$
Substituting these relations into \eqref{eq:adjoint-verlet} yields
\[
\boldsymbol{\lambda}_p^{\,n+\frac12}
=
\boldsymbol{\lambda}_p^{\,n+1},
\qquad
\boldsymbol{\lambda}_v^{\,n+\frac12}
=
\boldsymbol{\lambda}_v^{\,n+1}
+
\frac{\Delta t}{2}
\boldsymbol{\lambda}_p^{\,n+1},
\]
together with
\[
\mathbf M_{\mathrm L}
\frac{
\boldsymbol{\lambda}_p^{\,n+\frac12}
-
\boldsymbol{\lambda}_p^{\,n}
}
{\Delta t}
=
\mathbf K
\boldsymbol{\lambda}_v^{\,n+\frac12},
\qquad
\mathbf M_{\mathrm L}
\frac{
\boldsymbol{\lambda}_v^{\,n}
-
\boldsymbol{\lambda}_v^{\,n+\frac12}
}
{\Delta t}
=
-
\mathbf B
\boldsymbol{\lambda}_v^{\,n+\frac12}
+
\frac{1}{2}
\mathbf M_{\mathrm L}
\boldsymbol{\lambda}_p^{\,n}.
\]
These equations possess the standard mixed finite element structure
%\[
%\mathbf M_{\mathrm L}\dot{\boldsymbol{\lambda}}
%+\mathbf K\boldsymbol{\lambda}
%=[{\text{boundary terms}}],
%\]
which motivates the continuous interpretation of the algebraic transpose as a reverse-time first-order wave equation. The only remaining discrepancy with the continuous adjoint is the boundary contribution involving $\mathbf B$, which arises solely from the artificial truncation of the computational domain.
% ============================================================

% Your manuscript starts here.

% ============================================================
% Acknowledgments
% ============================================================

\section*{Acknowledgments}

B. M. Afkham and A. M. A. Alghamdi are partially funded by Villum Investigator grant no.\ 25893. B. M. Afkham is also supported by the Research Council of Finland under grant number 371523. H. Yazdanian is supported by Flagship of Advanced Mathematics for Sensing, Imaging and Modelling grant no.\ 2402070 and the Finnish Ministry of Education and Culture through the Doctoral Education Pilot for Mathematics of Sensing, Imaging and Modelling. T. Tarvainen is supported by the European Research Council grant agreement no.\ 101001417, and by the Research Council of Finland grant no.\ 353086, Flagship of Advanced Mathematics for Sensing Imaging and Modelling grant no.\ 358944, and Flagship Program Photonics Research and Innovation grant no.\ 320166.

% ============================================================
% AI declaration
% ============================================================

\section*{Use of Generative-AI tools declaration}

The authors used generative AI tools solely for language editing and improvement of the presentation of the manuscript. The scientific ideas, methodology, analysis, interpretations, and results presented in this work are entirely those of the authors. The authors reviewed and approved all AI-assisted edits and take full responsibility for the content of the manuscript.

% ============================================================
% References
% ============================================================

\bibliographystyle{plain}
\bibliography{references}

@book{rasmussen2006,
	author={C.~E. Rasmussen and C.~ K.~ I.~ Williams},
	title={Gaussian Processes for Machine Learning},
	publisher={MIT Press},
	year={2006}
}

@article{cox2005fast,
  title={Fast calculation of pulsed photoacoustic fields in fluids using k-space methods},
  author={Cox, Ben T and Beard, Paul C},
  journal={The Journal of the Acoustical Society of America},
  volume={117},
  number={6},
  pages={3616--3627},
  year={2005},
  publisher={Acoustical Society of America}
}

@Book{kaipio05,
  author = 	 {J. Kaipio and E. Somersalo},
  ALTeditor = 	 {},
  title = 	 {Statistical and Computational Inverse Problems},
  publisher = 	 {Springer},
  year = 	 {2005},
  OPTkey = 	 {},
  OPTvolume = 	 {},
  OPTnumber = 	 {},
  OPTseries = 	 {},
  address = 	 {New York},
  OPTedition = 	 {},
  OPTmonth = 	 {},
  OPTnote = 	 {},
  OPTannote = 	 {}
}

@article{sahlstrom2023utilizing,
  title={Utilizing variational autoencoders in the Bayesian inverse problem of photoacoustic tomography},
  author={Sahlstr{\"o}m, Teemu and Tarvainen, Tanja},
  journal={SIAM Journal on Imaging Sciences},
  volume={16},
  number={1},
  pages={89--110},
  year={2023},
  publisher={SIAM}
}

@Article{Goh2021,
  title =	 {Solving Bayesian Inverse Problems via Variational Autoencoders},
  author =	 {H. Goh and S. Sheriffdeen and J. Wittmer and T. Bui-Thanh},
  doi =		 {},
  journal =	 {Proceedings of Machine Learning Research},
  number =	 {},
  volume =	 {145},
  year =	 2021,
  month =	 {},
  pages =	 {386--425},
}

@article{tick2019modelling,
  title={Modelling of errors due to speed of sound variations in photoacoustic tomography using a Bayesian framework},
  author={Tick, Jenni and Pulkkinen, Aki and Tarvainen, Tanja},
  journal={Biomedical physics \& engineering express},
  volume={6},
  number={1},
  pages={015003},
  year={2019},
  publisher={IOP Publishing}
}

@article{hauptmann2020deep,
  title={Deep learning in photoacoustic tomography: current approaches and future directions},
  author={Hauptmann, Andreas and Cox, Ben},
  journal={Journal of Biomedical Optics},
  volume={25},
  number={11},
  pages={112903--112903},
  year={2020}
}

@Article{tick2016,
  author = 	 {J. Tick and A. Pulkkinen and T. Tarvainen},
  title = 	 {Image reconstruction with uncertainty quantification in photoacoustic tomography},
  journal =  jasa,
  year = 	 {2016},
  OPTkey = 	 {},
  volume = 	 {139},
  OPTnumber = 	 {},
  pages = 	 {1951-1961},
  OPTmonth = 	 {},
  OPTnote = 	 {},
  OPTannote = 	 {}
}

@Article{pulkkinen2014,
  author = 	 {A. Pulkkinen and B.~T. Cox and S.~R. Arridge and J.~P. Kaipio and T. Tarvainen},
  title = 	 {A {B}ayesian approach to spectral quantitative photoacoustic tomography},
  journal = 	 ip,
  year = 	 {2014},
  OPTkey = 	 {},
  volume = 	 {30},
  OPTnumber = 	 {},
  pages = 	 {065012},
  OPTmonth = 	 {},
  OPTnote = 	 {},
  OPTannote = 	 {}
}

@Article{beard2011,
  author = 	 {P. Beard},
  title = 	 {Biomedical photoacoustic imaging},
  journal = 	 {Interface Focus},
  year = 	 {2011},
  OPTkey = 	 {},
  volume = 	 {1},
  OPTnumber = 	 {},
  pages = 	 {602-631},
  OPTmonth = 	 {},
  OPTnote = 	 {},
  OPTannote = 	 {}
}

@Article{li2009,
  author = 	 {C. Li and L.~V. Wang},
  title = 	 {Photoacoustic tomography and sensing in biomedicine},
  journal = 	 pmb,
  year = 	 {2009},
  OPTkey = 	 {},
  volume = 	 {54},
  OPTnumber = 	 {},
  pages = 	 {R59-R97},
  OPTmonth = 	 {},
  OPTnote = 	 {},
  OPTannote = 	 {}
}

@Misc{amsmath,
  author =	 {{American Mathematical Society}},
  title =	 {User's Guide for the \texttt{amsmath} Package
                  (Version 2.0)},
  url =		 {ftp://ftp.ams.org/pub/tex/doc/amsmath/amsldoc.pdf},
  urldate =	 {2015-07-30},
  year =	 2002}

@article{hairer2006geometric,
  title={Geometric numerical integration},
  author={Hairer, Ernst and Hochbruck, Marlis and Iserles, Arieh and Lubich, Christian},
  journal={Oberwolfach Reports},
  volume={3},
  number={1},
  pages={805--882},
  year={2006}
}

@article{bardsley2014randomize,
  title={Randomize-then-optimize: A method for sampling from posterior distributions in nonlinear inverse problems},
  author={Bardsley, Johnathan M and Solonen, Antti and Haario, Heikki and Laine, Marko},
  journal={SIAM Journal on Scientific Computing},
  volume={36},
  number={4},
  pages={A1895--A1910},
  year={2014},
  publisher={SIAM}
}

@book{hansen1998rank,
  title={Rank-deficient and discrete ill-posed problems: numerical aspects of linear inversion},
  author={Hansen, Per Christian},
  year={1998},
  publisher={SIAM}
}

@article{xia2014photoacoustic,
  title={Photoacoustic tomography: principles and advances},
  author={Xia, Jun and Yao, Junjie and Wang, Lihong V},
  journal={Electromagnetic waves (Cambridge, Mass.)},
  volume={147},
  pages={1},
  year={2014},
  publisher={NIH Public Access}
}

@article{tarvainen2024quantitative,
  title={Quantitative photoacoustic tomography: modeling and inverse problems},
  author={Tarvainen, Tanja and Cox, Ben},
  journal={Journal of Biomedical Optics},
  volume={29},
  number={S1},
  pages={S11509--S11509},
  year={2024},
  publisher={Society of Photo-Optical Instrumentation Engineers}
}

@article{tarvainen2012reconstructing,
  title={Reconstructing absorption and scattering distributions in quantitative photoacoustic tomography},
  author={Tarvainen, Tanja and Cox, Benjamin T and Kaipio, JP and Arridge, Simon R},
  journal={Inverse Problems},
  volume={28},
  number={8},
  pages={084009},
  year={2012},
  publisher={IOP Publishing}
}

@article{tarvainen2013bayesian,
  title={Bayesian image reconstruction in quantitative photoacoustic tomography},
  author={Tarvainen, Tanja and Pulkkinen, Aki and Cox, Ben T and Kaipio, Jari P and Arridge, Simon R},
  journal={IEEE transactions on medical imaging},
  volume={32},
  number={12},
  pages={2287--2298},
  year={2013},
  publisher={IEEE}
}

@article{doi:10.1137/18M1231341,
author = {Ren, Kui and Vall\'{e}lian, Sarah},
title = {Characterizing Impacts of Model Uncertainties in Quantitative Photoacoustics},
journal = {SIAM/ASA Journal on Uncertainty Quantification},
volume = {8},
number = {2},
pages = {636-667},
year = {2020},
doi = {10.1137/18M1231341},
URL = {https://doi.org/10.1137/18M1231341},
eprint = {https://doi.org/10.1137/18M1231341}
}

@article{sahlstrom2020modeling,
  title={Modeling of errors due to uncertainties in ultrasound sensor locations in photoacoustic tomography},
  author={Sahlstr{\"o}m, Teemu and Pulkkinen, Aki and Tick, Jenni and Leskinen, Jarkko and Tarvainen, Tanja},
  journal={IEEE Transactions on Medical Imaging},
  volume={39},
  number={6},
  pages={2140--2150},
  year={2020},
  publisher={IEEE}
}

@article{roininen2014whittle,
  title={WHITTLE-MAT{\'E}RN PRIORS FOR BAYESIAN STATISTICAL INVERSION WITH APPLICATIONS IN ELECTRICAL IMPEDANCE TOMOGRAPHY.},
  author={Roininen, Lassi and Huttunen, Janne MJ and Lasanen, Sari},
  journal={Inverse Problems \& Imaging},
  volume={8},
  number={2},
  year={2014}
}

@article{engquist1979radiation,
  title={Radiation boundary conditions for acoustic and elastic wave calculations},
  author={Engquist, Bjorn and Majda, Andrew},
  journal={Communications on pure and applied mathematics},
  volume={32},
  pages={313--357},
  year={1979}
}

@article{wang2017bayesian,
  title={Bayesian inverse problems with $l\_1$ priors: a randomize-then-optimize approach},
  author={Wang, Zheng and Bardsley, Johnathan M and Solonen, Antti and Cui, Tiangang and Marzouk, Youssef M},
  journal={SIAM Journal on Scientific Computing},
  volume={39},
  number={5},
  pages={S140--S166},
  year={2017},
  publisher={SIAM}
}

@article{paige1982lsqr,
  title={LSQR: An algorithm for sparse linear equations and sparse least squares},
  author={Paige, Christopher C and Saunders, Michael A},
  journal={ACM Transactions on Mathematical Software (TOMS)},
  volume={8},
  number={1},
  pages={43--71},
  year={1982},
  publisher={ACM New York, NY, USA}
}

@inproceedings{seoni2025exploring,
  title={Exploring uncertainty quantification for photoacoustic image reconstruction and quantitative oxygenation mapping},
  author={Seoni, Silvia and Scardigno, Roberto M and Cotrufo, Bruna and Salvi, Massimo and Brunetti, Antonio and Guerriero, Andrea and Rotunno, Giulia and Buongiorno, Domenico and Vallan, Alberto and Molinari, Filippo and others},
  booktitle={Photons Plus Ultrasound: Imaging and Sensing 2025},
  volume={13319},
  pages={265--270},
  year={2025},
  organization={SPIE}
}

@article{dong2019fixing,
  title={Fixing nonconvergence of algebraic iterative reconstruction with an unmatched backprojector},
  author={Dong, Yiqiu and Hansen, Per Christian and Hochstenbach, Michiel E and Brogaard Riis, Nicolai André},
  journal={SIAM Journal on Scientific Computing},
  volume={41},
  number={3},
  pages={A1822--A1839},
  year={2019},
  publisher={SIAM}
}

@inproceedings{huber2025novel,
  title={A Novel Interpretation of the Radon Transform’s Ray and Pixel-Driven Discretizations Under Balanced Resolutions},
  author={Huber, Richard},
  booktitle={International Conference on Scale Space and Variational Methods in Computer Vision},
  pages={132--145},
  year={2025},
  organization={Springer}
}

@article{fried1975finite,
  title={Finite element mass matrix lumping by numerical integration with no convergence rate loss},
  author={Fried, Isaac and Malkus, David S},
  journal={International Journal of Solids and Structures},
  volume={11},
  number={4},
  pages={461--466},
  year={1975},
  publisher={Elsevier}
}

@misc{afkham2026pat3duq,
  author       = {Afkham, Babak Maboudi and
                  Alghamdi, Amal Mohammed A. and
                  Yazdanian, Hassan and
                  Tarvainen, Tanja},
  title        = {{PAT-3D-Uncertainty-Quantification}: Accompanying Code for
                  3D Uncertainty Quantification for Photoacoustic Tomography},
  year         = {2026},
  howpublished = {\url{https://github.com/babakmaboudi/PAT-3D-Uncertainty-Quantification}},
  note         = {GitHub repository}
}

\end{document}